\documentclass[11pt, reqno]{amsart}

\usepackage[usenames,dvipsnames]{color}

\usepackage[OT2, T1]{fontenc}
\usepackage{url}
\usepackage{amsmath}
\usepackage{array}
\usepackage{graphicx}
\usepackage{amssymb}
\usepackage{amsthm}
\definecolor{link}{RGB}{11,0,128}
\usepackage[colorlinks=true,citecolor=NavyBlue,linkcolor=OliveGreen,urlcolor=link]{hyperref}
\usepackage{hyperref}
\usepackage{paralist}
\usepackage{colonequals}			
\usepackage{color}
\usepackage{mathabx}			
\usepackage{amscd}
\usepackage[all,cmtip]{xy}			
\usepackage{sseq}				
\usepackage{verbatim}
\usepackage{parskip}
\usepackage{microtype}
\usepackage[in]{fullpage}
\usepackage{mathrsfs} 			
\usepackage{cleveref}
\usepackage{enumitem}
\usepackage{moreenum}			
\usepackage[alphabetic,lite]{amsrefs} 	
\usepackage{stmaryrd}			
\usepackage[foot]{amsaddr}
\usepackage{subfiles}
\usepackage{ragged2e}
\usepackage{stackengine}

\newcommand{\gA}{\alpha}

\newcommand{\bA}{\mathbb{A}}

\newcommand{\bF}{\mathbb{F}}
\newcommand{\bG}{\mathbb{G}}

\newcommand{\bP}{\mathbb{P}}
\newcommand{\bQ}{\mathbb{Q}}

\newcommand{\bZ}{\mathbb{Z}}

\newcommand{\cE}{\mathcal{E}}
\newcommand{\cF}{\mathcal{F}}
\newcommand{\cG}{\mathcal{G}}

\newcommand{\cO}{\mathcal{O}}

\newcommand{\cT}{\mathcal{T}}

\newcommand{\cY}{\mathcal{Y}}

\newcommand{\fm}{\mathfrak{m}}

\newcommand{\fp}{\mathfrak{p}}

\newcommand{\sB}{\mathscr{B}}

\newcommand{\sE}{\mathscr{E}}
\newcommand{\sF}{\mathscr{F}}
\newcommand{\sG}{\mathscr{G}}

\newcommand{\sO}{\mathscr{O}}
\newcommand{\sP}{\mathscr{P}}

\newcommand{\sR}{\mathscr{R}}

\newcommand{\sT}{\mathscr{T}}

\newcommand{\sY}{\mathscr{Y}}
\newcommand{\sZ}{\mathscr{Z}}

\DeclareMathOperator{\Bl}{Bl}			
\DeclareMathOperator{\Frac}{Frac}		
\DeclareMathOperator{\GL}{GL}		
\DeclareMathOperator{\Isom}{Isom}		
\DeclareMathOperator{\Ker}{Ker}		
\DeclareMathOperator{\PGL}{PGL}		
\DeclareMathOperator{\Pic}{Pic}		
\DeclareMathOperator{\rad}{rad}		
\DeclareMathOperator{\Res}{Res}		
\DeclareSymbolFont{cyrletters}{OT2}{wncyr}{m}{n}
\DeclareMathSymbol{\Sha}{\mathalpha}{cyrletters}{"58}	
\DeclareMathOperator{\Spec}{Spec}		

\newcommand{\ad}{\mathrm{ad}}			

\newcommand{\ce}{\colonequals}

\newcommand{\eps}{\varepsilon}

\newcommand{\gp}{{\mathrm{gp}}}		
\newcommand{\hra}{\hookrightarrow}
\renewcommand{\i}{^{-1}}
\newcommand{\id}{\mathrm{id}}			

\newcommand{\isomto}{\overset{\sim}{\longrightarrow}}

\newcommand{\Nis}{\mathrm{Nis}}	
\newcommand{\ov}{\overline}

\newcommand{\ra}{\rightarrow}

\newcommand{\red}{{\mathrm{red}}}		

\newcommand{\surjects}{\twoheadrightarrow}
\newcommand{\tensor}{\otimes} 			

\newcommand{\wt}{\widetilde}
\newcommand{\xra}{\xrightarrow}
\newcommand{\Zar}{\mathrm{Zar}}		

\providecommand{\up}[1]{{\upshape(}#1{\upshape)}}
\providecommand{\uref}[1]{{\upshape\ref{#1}}}

\providecommand{\f}[2]{\frac{#1}{#2}}

\renewcommand{\b}{\textbf}
\providecommand{\ucolon}{{\upshape:} }
\providecommand{\uscolon}{{\upshape;} }

\newcommand{\brems}{\begin{rems} \hfill \begin{enumerate}[label=\b{\thenumberingbase.},ref=\thenumberingbase]}
\newcommand{\remi}{\addtocounter{numberingbase}{1} \item}
\newcommand{\erems}{\end{enumerate} \end{rems}}
\newcommand{\begs}{\begin{egs} \hfill \begin{enumerate}[label=\b{\thenumberingbase.},ref=\thenumberingbase]}

\newcommand{\eegs}{\end{enumerate} \end{egs}}
\newcommand{\m}{\item}
\newcommand{\bsm}{\begin{smallmatrix}}
\newcommand{\esm}{\end{smallmatrix}}
\newcommand{\blem}{\begin{lemma}}
\newcommand{\elem}{\end{lemma}}

\newcommand{\bconj}{\begin{conj}}
\newcommand{\econj}{\end{conj}}
\newcommand{\bprob}{\begin{Problem}}
\newcommand{\eprob}{\end{Problem}}

\newcommand{\bq}{\begin{Q}}
\newcommand{\eq}{\end{Q}}
\newcommand{\benum}{\begin{enumerate}[label={{\upshape(\alph*)}}]}
\newcommand{\benuma}{\begin{enumerate}[label={{\upshape(\arabic*)}}]}
\newcommand{\benumr}{\begin{enumerate}[label={{\upshape(\roman*)}}]}
\newcommand{\eenum}{\end{enumerate}}
\newcommand{\bitem}{\begin{itemize}}
\newcommand{\eitem}{\end{itemize}}
\newcommand{\bc}{}
\newcommand{\bd}{\begin{defn}}
\newcommand{\ed}{\end{defn}}

\newcommand{\beg}{\begin{eg}}
\newcommand{\eeg}{\end{eg}}

\newcommand{\bcl}{\begin{claim}}
\newcommand{\ecl}{\end{claim}}

\newcommand{\x}{\text}

\newcommand{\q}{\quad}
\providecommand{\qxq}[1]{\quad\text{#1}\quad}
\providecommand{\qx}[1]{\quad\text{#1}}

\newcommand{\qq}{\quad\quad}
\newcommand{\qqq}{\quad\quad\quad}

\newcommand{\tst}{\textstyle}
\newcommand{\ba}{\begin{aligned}}
\newcommand{\ea}{\end{aligned}}
\newcommand{\be}{\begin{equation}}
\newcommand{\ee}{\end{equation}}
\newcommand{\bpf}{\begin{proof}}
\newcommand{\epf}{\end{proof}}
\newcommand{\bthm}{\begin{thm}}
\newcommand{\ethm}{\end{thm}}
\newcommand{\bprop}{\begin{prop}}
\newcommand{\eprop}{\end{prop}}
\newcommand{\bcor}{\begin{cor}}
\newcommand{\ecor}{\end{cor}}
\newcommand{\brem}{\begin{rem}}
\newcommand{\erem}{\end{rem}}
\newcommand*{\QED}{\hfill\ensuremath{\qed}}

\usepackage{aliascnt}
\newaliascnt{numberingbase}{subsection}
\numberwithin{equation}{numberingbase}

\newtheoremstyle{thms}{0.5em}{0pt}{\itshape}{}{\bfseries}{.}{ }{}
\theoremstyle{thms}
\newtheorem{conj}[numberingbase]{Conjecture}

\newtheorem{cor}[numberingbase]{Corollary}

\newtheorem{lemma}[numberingbase]{Lemma}

\newtheorem{prop}[numberingbase]{Proposition}

\newtheorem{Q}[numberingbase]{Question}
\newtheorem{thm}[numberingbase]{Theorem}

\newtheoremstyle{claims}{0.5em}{0pt}{}{}{\itshape}{.}{ }{}
\theoremstyle{claims}
\newtheorem{claim}[equation]{Claim}

\newtheoremstyle{defs}{0.5em}{0pt}{}{}{\bfseries}{.}{ }{}
\theoremstyle{defs}
\newtheorem{defn}[numberingbase]{Definition}

\newtheorem{eg}[numberingbase]{Example}
\newtheorem*{egs}{Examples}
\newtheorem{rem}[numberingbase]{Remark}

\newtheorem*{rems}{Remarks}

\Crefname{claim}{Claim}{Claims}
\Crefname{bclaim}{Claim}{Claims}
\Crefname{sublemma}{Lemma}{Lemmas}
\Crefname{conj}{Conjecture}{Conjectures}
\Crefname{cor}{Corollary}{Corollaries}
\Crefname{defn}{Definition}{Definitions}
\Crefname{eg}{Example}{Examples}
\Crefname{prop}{Proposition}{Propositions} 
\Crefname{Q}{Question}{Questions}
\Crefname{rem}{Remark}{Remarks}
\Crefname{thm}{Theorem}{Theorems}
\Crefname{Theorem}{Theorem}{Theorems}
\Crefname{variant}{Variant}{Variants}
\Crefname{caution}{Caution}{Cautions}

\theoremstyle{thms}
\newtheorem{thm-tweak}[subsection]{Theorem}
\Crefname{thm-tweak}{Theorem}{Theorems}
\newtheorem{lemma-tweak}[subsection]{Lemma}
\Crefname{lemma-tweak}{Lemma}{Lemmas}
\newtheorem{cor-tweak}[subsection]{Corollary}
\Crefname{cor-tweak}{Corollary}{Corollaries}
\newtheorem{prop-tweak}[subsection]{Proposition}
\Crefname{prop-tweak}{Proposition}{Propositions} 
\newtheorem{conj-tweak}[subsection]{Conjecture}
\Crefname{conj-tweak}{Conjecture}{Conjectures} 
\newtheorem{q-tweak}[subsection]{Question}
\Crefname{q-tweak}{Question}{Questions} 

\theoremstyle{defs}
\newtheorem{defn-tweak}[subsection]{Definition}
\Crefname{defn-tweak}{Definition}{Definitions}
\newtheorem{eg-tweak}[subsection]{Example}
\Crefname{eg-tweak}{Example}{Examples}
\newtheorem*{rems-tweak}{Remarks}
\newtheorem{rem-tweak}[subsection]{Remark}
\Crefname{rem-tweak}{Remark}{Remarks}

\newtheoremstyle{subsection-tweak}
   {2pt}
   {3pt}%
   {}
   {}%
   {\bfseries}
   {}%
   {.5em}
   {\thmnumber{\@{#1}{}\@{#2}.}%
    \thmnote{~{\bfseries#3.}}}    
    
\theoremstyle{subsection-tweak}
\newtheorem{pp}[numberingbase]{}
\newcommand{\bpp}{\begin{pp}}
\newcommand{\epp}{\end{pp}}

\theoremstyle{subsection-tweak}
\newtheorem{pp-tweak}[subsection]{}

\makeatletter
\def\@tocline#1#2#3#4#5#6#7{
    \begingroup 
    \@ifempty{#4}{}{}

    \parindent\z@ \leftskip#3\relax \advance\leftskip\@tempdima\relax
    #5\hskip-\@tempdima
      \ifcase #1
       \or\or \hskip 2em \or \hskip 1em \else \hskip 3em \fi%
      #6\nobreak\relax
    \dotfill\hbox to\@pnumwidth{\@tocpagenum{#7}}\par
    \nobreak
    \endgroup
 }
 \def\l@section{\@tocline{1}{0pt}{1pc}{}{}}

\renewcommand{\tocsection}[3]{%
  \indentlabel{\@ifnotempty{#2}{\makebox[1.3em][l]{%
    \ignorespaces#1 \bfseries{#2}.\hfill}}}\bfseries{#3}
    \vspace{-5pt}}

\renewcommand{\tocsubsection}[3]{%
  \indentlabel{\@ifnotempty{#2}{\hspace*{-0.5em}\makebox[2.1em][l]{%
    \ignorespaces#1#2.\hfill}}}#3
    \vspace{-5pt}}
   
\makeatother 

\makeatletter 
\newcommand\appendix@section[1]{%
  \refstepcounter{section}%
  \orig@section*{Appendix \@Alph\c@section. #1}%
}
\let\orig@section\section
\g@addto@macro\appendix{\let\section\appendix@section}
\makeatother

\makeatletter
\@namedef{subjclassname@2020}{%
  \textup{2020} Mathematics Subject Classification}
\makeatother

\author{K\k{e}stutis \v{C}esnavi\v{c}ius}
\address{CNRS, Universit\'{e} Paris-Saclay,   Laboratoire de math\'{e}matiques d'Orsay, F-91405, Orsay, France}
\email{kestutis.cesnavicius@cnrs.fr}

\date{\today}

\begin{document}


\title{Torsors on the complement of a smooth divisor}

\maketitle

\begin{abstract}
We complete the proof of the Nisnevich conjecture in equal characteristic: for a smooth algebraic variety $X$ over a field $k$, a $k$-smooth divisor $D \subset X$, and a reductive $X$-group $G$ whose base change $G_D$ is totally isotropic, we show that each generically trivial $G$-torsor on $X\setminus D$ trivializes Zariski semilocally on $X$. In mixed characteristic, we show the same when $k$ is a replaced by a discrete valuation ring $\cO$, the divisor $D$ is the closed $\cO$-fiber of $X$, and either $G$ is quasi-split or $G$ is only defined over $X \setminus D$ but descends to a quasi-split group over $\Frac(\cO)$ (a Kisin--Pappas type variant). 
Our arguments combine Gabber--Quillen style presentation lemmas with excision and reembedding d\'{e}vissages to reduce to analyzing generically trivial torsors over a relative affine line. As a byproduct of this analysis, we give a new proof for the Bass--Quillen conjecture for reductive group torsors over $\bA^d_R$ in equal characteristic.

\end{abstract}

\hypersetup{
    linktoc=page,     
}
\renewcommand*\contentsname{}
\tableofcontents

\section{The corrected statement of the Nisnevich conjecture and our main results}

In \cite{Nis89}*{Conjecture 1.3}, Nisnevich proposed a common generalization of the Quillen conjecture \cite{Qui76}*{(2) on page 170} that had grown out of Serre's problem about vector bundles on affine spaces and of the Grothendieck--Serre conjecture \cite{Ser58b}*{page 31, Remarque}, \cite{Gro58}*{pages 26--27, Remarques 3} about Zariski local triviality of generically trivial torsors under reductive groups. In its geometric case, the Nisnevich conjecture predicts that, for a reductive group scheme $G$ over a smooth variety $X$ over a field $k$ and a $k$-smooth divisor $D \subset X$, every generically trivial $G$-torsor on $X \setminus D$ trivializes Zariski locally on $X$. Recent counterexamples of Fedorov \cite{Fed22b}*{Proposition~4.1} show that this fails for anisotropic $G$, so, to bypass them,  one considers the following isotropicity condition whose relevance for problems about torsors has been observed already in \cite{Rag89}.

\bd[\cite{split-unramified}*{Definition 8.1}] \label{def:totally-isotropic}
Let $S$ be a scheme and let $G$ be a reductive $S$-group scheme. We say that $G$ is \emph{totally isotropic}\footnote{In \cite{Fed22a} and \cite{Fed22b}, the terminology `strongly locally isotropic' was used for the same notion.} if in the canonical~decomposition 
\be \label{eqn:Gad-decomposition}
\tst G^\ad \cong  \prod_{i\in \{A_n,\, B_n,\, \dotsc,\, G_2\}} \Res_{S_i/S}(G_i)
\ee
of \cite{SGA3IIInew}*{expos\'{e} XXIV, proposition 5.10 (i)}, in which $i$ ranges over the types of connected Dynkin diagrams, $S_i$ is a finite \'{e}tale $S$-scheme, and $G_i$ is an adjoint semisimple $S_i$-group with simple geometric fibers of type $i$, Zariski locally on $S$ each $G_i$ has $\bG_{m,\, S_i}$ as a subgroup. 
\ed

Intuitively, $G$ is totally isotropic if and only if its simple factors are isotropic. Recall from \cite{SGA3IIInew}*{expos\'{e}~XXVI, corollaire~6.12} that  in \Cref{def:totally-isotropic} it is equivalent to require that Zariski locally on $S$ each $G_i$ has a parabolic $S_i$-subgroup that contains no $S_i$-fiber of $G_i$. For instance, every quasi-split, so also every split, group is totally isotropic, as is every torus. 

With the total isotropicity in place, the Nisnevich conjecture becomes the following statement.

\bconj[Nisnevich] \label{conj:Nisnevich}
For a  regular semilocal ring $R$, an $r \in R$ that is a regular parameter in the sense that $r \not \in \fm^2$ for each maximal ideal $\fm \subset R$, and a reductive $R$-group scheme $G$ such that $G_{R/(r)}$ is totally isotropic, every generically trivial $G$-torsor over $R[\f1r]$ is trivial, that is,
\[
\tst \Ker(H^1(R[\f1r], G) \ra H^1(\Frac(R), G)) = \{*\}.
\]
\econj

For instance, in the case when $r$ is a unit, the total isotropicity condition holds for every reductive $R$-group $G$ and we recover the Grothendieck--Serre conjecture. The condition also holds in the case when $G$ is a torus, and this case follows from the known toral case of the Grothendieck--Serre conjecture, see \cite{torsors-regular}*{Section 3.4.2~(1)}. In \cite{Fed22b}, Fedorov settled the Nisnevich conjecture in the case when $R$ contains an infinite field and $G$ itself is totally isotropic. Other than this, some low dimensional cases are known, see \cite{torsors-regular}*{Section~3.4.2}---for instance, the case when $R$ is local of dimension $\le 3$ and $G$ is either $\GL_n$ or $\PGL_n$ is a result of Gabber \cite{Gab81}*{Chapter I, Theorem 1}. 

We settle the Nisnevich conjecture in equal characteristic and in some mixed characteristic cases.

\bthm \label{thm:main-DVR}
Let $R$ be a regular semilocal ring, let $r \in R$ be a regular parameter in the sense that $r \not\in \fm^2$ for each maximal ideal $\fm \subset R$, and let $G$ be a reductive $R[\f1r]$-group. 
In the following cases, 
\[
\tst \Ker(H^1(R[\f1r], G) \ra H^1(\Frac(R), G)) = \{*\},
\]
in other words, in the following cases every generically trivial $G$-torsor over $R[\f1r]$ is trivial\ucolon
\benuma
\m \label{m:mDVR-1}
{\upshape(\S\ref{pp:Nisnevich-conj-pf})}
if $R$ contains a field and $G$ extends to a reductive $R$-group $\cG$ with $\cG_{R/(r)}$ totally~isotropic\uscolon


\m \label{m:mDVR-3}
{\upshape(\S\ref{pp:pf-mixed-char})}
if $R$ is geometrically regular\footnote{For a ring $A$, recall that an $A$-algebra $B$ is  \emph{geometrically regular} if it is flat and the base change of each of its $A$-fibers to any finite field extension of the corresponding residue field of $A$ is regular, see \cite{SP}*{Definition~\href{https://stacks.math.columbia.edu/tag/0382}{0382}}. For instance, $R$ could be a semilocal ring of a smooth algebra over a discrete valuation ring $\cO$ with $r$ as a uniformizer.} over a Dedekind subring $\cO$ containing $r$ and $G$ either extends to a quasi-split reductive $R$-group or descends to a quasi-split reductive $\cO[\f1r]$-group.
\eenum
\ethm

The mixed characteristic case \ref{m:mDVR-3} is new already for vector bundles, that is, for $G = \GL_n$. In contrast, at least for local $R$, the vector bundle case of the equicharacteristic \ref{m:mDVR-1} is due to Bhatwadekar--Rao \cite{BR83}*{Theorem 2.5}. When $r \in R^\times$, \Cref{thm:main-DVR} recovers the equal and mixed characteristic cases of the Grothendieck--Serre conjecture settled in \cite{FP15}, \cite{Pan20a}, \cite{split-unramified}, so we reprove these here.

In the mixed characteristic case \ref{m:mDVR-3}, the requirement that $r \in \cO$ is quite restrictive relative to the assumptions of \Cref{conj:Nisnevich}. However, the case of \ref{m:mDVR-3} in which $G$ descends to an $\cO[\f1r]$-group but need not extend to a reductive $R$-group was inspired by Kisin--Pappas \cite{KP18}*{Section~1.4, especially, Lemma 1.4.6}, who used such a statement for some $2$-dimensional $R$ under further assumptions on $G$. 

The geometric version of \Cref{thm:main-DVR}~\ref{m:mDVR-1} is the following statement announced in the abstract. 

\bthm\label{thm:main-G}
For a field $k$, a smooth $k$-scheme $X$, a $k$-smooth divisor $D \subset X$, and a reductive $X$-group scheme $G$ such that $G_D$ is totally isotropic, every generically trivial $G$-torsor $E$ over $X \setminus D$ is trivial Zariski semilocally on $X$, that is, for every $x_1, \dotsc, x_m \in X$ that lie in a single affine open, there is an affine open $U \subset X$ containing all the $x_i$ such that $E|_{U \setminus D}$ is trivial. 
\ethm

\Cref{thm:main-G} follows by applying \Cref{thm:main-DVR}~\ref{m:mDVR-1} to the semilocal ring of $X$ at $x_1, \dotsc, x_m$ (built via prime avoidance, see \cite{SP}*{Lemma~\href{https://stacks.math.columbia.edu/tag/00DS}{00DS}}) and spreading out. 
Even when $X$ is affine,  the stronger statement that $E$ extends to a $G$-torsor over $X$ is false: for $G = \GL_n$, this had been a question of Quillen \cite{Qui76}*{(3) on page 170} that was answered negatively by Swan in \cite{Swa78}*{Section~2}. Even for $\GL_n$, \Cref{thm:main-G} typically fails if $D$ is singular or if $X$ is singular, see \cite{Lam06}*{pages 34--35}.

We use \Cref{thm:main-DVR} to reprove the following equal characteristic case of the generalization of the Bass--Quillen conjecture to torsors under reductive group schemes \cite{torsors-regular}*{Conjecture 3.6.1}.




\bthm[\S\ref{pp:BQ-pf}] \label{thm:BQ-ann}
For a regular ring $R$ containing a field and a totally isotropic 
reductive $R$-group scheme $G$, every generically trivial $G$-torsor over $\bA^d_R$ descends to a $G$-torsor over $R$, equivalently,
\[
H^1_\Zar(R, G) \isomto H^1_{\Zar}(\bA^d_R, G) \qxq{or, if one prefers,} H^1_\Nis(R, G) \isomto H^1_{\Nis}(\bA^d_R, G).
\]
\ethm

The equivalence of the three formulations in \Cref{thm:BQ-ann} follows from the Grothendieck--Serre conjecture: by \Cref{thm:main-DVR}, a $G$-torsor over $\bA^d_R$ is generically trivial, if and only if it is Zariski locally trivial, if and only if it is Nisnevich locally trivial. The generic triviality assumption is needed because, for instance, for every separably closed field $k$ that is not algebraically closed, there are nontrivial $\PGL_n$-torsors over $\bA^1_k$, see \cite{CTS21}*{Theorem~5.6.1~(vi)}. The total isotropicity assumption is needed because of \cite{BS17a}*{Proposition~4.9}, where Balwe and Sawant show that a Bass--Quillen statement cannot hold beyond totally isotropic $G$. For earlier counterexamples to generalizations of the Bass--Quillen conjecture beyond totally isotropic reductive groups, see \cite{OS71}*{Propositions 1 and~2}, 
\cite{Par78} and \cite{Fed16}*{Theorem 3 (ii) (whose assumptions can be met by Remark 2.6~(i))}.


\Cref{thm:BQ-ann} was established by Stavrova in \cite{Sta22}*{Corollary 5.5} by a different method, and in the case when $R$ contains an infinite field already in the earlier \cite{Sta19}*{Theorem 4.4}. Prior to that, the case when $R$ is smooth over a field $k$ and $G$ is defined and totally isotropic over $k$ was settled by Asok--Hoyois--Wendt: they used methods of $\bA^1$-homotopy theory of Morel--Voevodsky to verify axioms of Colliot-Th\'{e}l\`{e}ne--Ojanguren \cite{CTO92} that were known to imply the statement, see \cite{AHW18}*{Theorem~3.3.7} for infinite $k$ and \cite{AHW20}*{Theorem 2.4} for finite $k$. 
For regular $R$ of mixed characteristic, \Cref{thm:BQ-ann} is only known in sporadic cases, for instance, when $G$ is a torus, see \cite{CTS87}*{Lemma 2.4}, as well as \cite{torsors-regular}*{Section~3.6.4} for an overview.


We obtain \Cref{thm:main-DVR} by refining the Grothendieck--Serre strategies used in \cite{Fed22b} and \cite{split-unramified}. In fact, we 
establish the following version of Grothendieck--Serre valid over arbitrary base rings. 

\bthm[\Cref{rem:AB}] \label{thm:abstract-GS-ann}
For a reductive group $G$ over a ring $A$, every $G$-torsor over a smooth affine $A$-curve $C$ that is trivial away from some $A$-finite $Z \subset C$ trivializes Zariski semilocally on $C$. 
\ethm

\Cref{thm:abstract-GS-ann}, more precisely, its finer version given in \Cref{thm:abstract-GS}, is our ultimate source of triviality of torsors under reductive groups, and it generalizes \cite{Fed22a}*{Theorem~4}, as well as several earlier results in the literature.  Armed with it we quickly reprove the cases of the Grothendieck--Serre conjecture that have been settled in \cite{FP15}, \cite{Pan20a}, \cite{split-unramified}: more precisely, we use Popescu approximation and presentation lemmas in the style of Gabber--Quillen to reduce these cases to the relative curve setting of \Cref{thm:abstract-GS-ann}, and in this way we dissect the overall argument into a part that works over arbitrary rings and a part that is specific to regular rings.

Coming back to the Nisnevich conjecture itself, a key novelty of our approach is the following extension result for $G$-torsors over smooth relative curves.


\bthm[Proposition \uref{prop:extend-0} and Theorem \uref{thm:rel-GS}] \label{thm:extend-Y-ann} Let $R$ be a regular semilocal ring containing a field and let $G$ be a reductive $R$-group. For a smooth affine $R$-scheme $C$ of pure relative dimension $1$ and an $R$-\up{finite \'{e}tale} closed $Y \subset C$ such that $G_Y$ is totally isotropic, every $G$-torsor $E$ over $C \setminus Y$ that is trivial away from some $R$-finite closed $Z \subset C$ extends to a $G$-torsor over $C$.
\ethm

Roughly, extending a $G$-torsor to all of $C$ in \Cref{thm:extend-Y-ann} corresponds to extending a $G$-torsor in \Cref{thm:main-DVR}~\ref{m:mDVR-1} to all of $R$, in effect, to reducing the Nisnevich conjecture to the Grothendieck--Serre conjecture---this is why \Cref{thm:extend-Y-ann} is crucial for us. Conversely, to reduce \Cref{thm:main-DVR}~\ref{m:mDVR-1} to \Cref{thm:extend-Y-ann} we use a presentation lemma that extends its variants due to Quillen and Gabber: we first use Popescu theorem to pass to the geometric setting of \Cref{thm:main-G} and then show in \Cref{lem:Gabber-Quillen} that, up to replacing $X$ by an affine open neighborhood of $x_1, \dotsc, x_m$, we can express $X$ as a smooth relative curve over some affine open of $\bA^{d - 1}_k$ in such a way that $D$ is relatively finite \'{e}tale and our generically trivial $G$-torsor over $X$ is trivial away from a relatively finite closed subscheme. 

As for \Cref{thm:extend-Y-ann}, in \S\ref{sec:extend-Y} we present a series of excision and patching d\'{e}vissages to reduce to when $C = \bA^1_R$ and $C \setminus Y$ descends to a smooth curve defined over a subfield $k\subset R$. In this ``constant'' case, we show that our $G$-torsor over $C\setminus Y$ is even trivial by the ``relative Grothendieck--Serre'' theorem of Fedorov from \cite{Fed22a} (with an earlier version due to Panin--Stavrova--Vavilov \cite{PSV15}) that we reprove in \Cref{thm:rel-GS}: for every $k$-algebra $W$, no nontrivial $G$-torsor over $R \tensor_k W$ trivializes over $\Frac(R) \tensor_k W$; the total isotropicity assumption is crucial for this beyond the ``classical'' case $W = \Spec(k)$. 
As for the excision and patching techniques, finite field obstructions are a well-known, delicate difficulty in the field. We overcome them with a novel version of Panin's ``finite field tricks'' presented in \Cref{lem:pump-Y}. The wide scope of these techniques makes our  overall approach to \Cref{thm:main-DVR} quite axiomatic, and although we do not pursue this here, it would be interesting to have similar results for other functors, for instance, for the unstable $K_1$-functor studied by Stavrova and her coauthors, compare, for instance, with  \cite{Sta22}, \cite{Sta19} and earlier articles cited there. 





\bpp[Notation and conventions] \label{conv}
All rings we consider are commutative and unital. For a point $s$ of a scheme (resp.,~for a prime ideal $\fp$ of a ring), we let $k_s$ (resp.,~$k_\fp$) denote its residue field. For a scheme $S$ over a ring $A$ and an $A$-algebra $B$, we write $S \tensor_A B$ for the base change $S \times_{\Spec A} \Spec B$. For a global section $s$ of a scheme $S$, we write $S[\f1s] \subset S$ for the open locus where $s$ does not vanish. For a ring $A$, we let $\Frac(A)$ denote its total ring of fractions. For a semilocal regular ring $R$, we say that an $r \in R$ is a \emph{regular parameter} if $r \not \in \fm^2$ for every maximal ideal $\fm \subset R$. 

For reductive groups, we use the terminology from SGA 3, as reviewed in \cite{torsors-regular}*{Section 1.3}. For a parabolic subgroup $P$ of a reductive group scheme $G$, we let $\sR_u(P)$ denote its unipotent radical constructed in \cite{SGA3IIInew}*{expos\'{e} XXVI, proposition 1.6 (i)}. We say that a torus $T$ over a scheme $S$ is \emph{isotrivial} if it splits over some finite \'{e}tale cover over $S$; this always holds if either $S$ is locally Noetherian and geometrically unibranch (in the sense that the map from the normalization of $S_\red$ to $S$ is a universal homeomorphism), see \cite{SGA3II}*{expos\'{e} X, th\'{e}or\`{e}me 5.16}, or if $T$ is of rank $\le 1$.
\epp

\subsection*{Acknowledgements} 
This article was inspired by the recent preprint \cite{Fed22b}, in which Roman Fedorov settled \Cref{thm:main-DVR}~\ref{m:mDVR-1} in the case when $R$ contains an infinite field and $G$ is totally isotropic. I thank him for a seminar talk on this subject and for helpful correspondence. I thank the referee for very helpful remarks that have significantly improved the manuscript. I thank Alexis Bouthier, Elden Elmanto, Ofer Gabber, Arnab Kundu, Shang Li, and Anastasia Stavrova for helpful conversations or correspondence, especially, Ofer Gabber for astute remarks during seminar talks. This project has received funding from the European Research Council (ERC) under the European Union's Horizon 2020 research and innovation programme (grant agreement No.~851146).


\section{Torsors over $\bA^d_A$} \label{sec:affine-Gr}



Our eventual source of triviality of torsors is the following general result about torsors over $\bP^1_A$. Its part \ref{m:CFR-b} is how the total isotropicity assumption ultimately enters into the geometric approach to the Nisnevich conjecture \ref{conj:Nisnevich} that is developed in this article building on \cite{Fed22b}. Earlier weaker versions of \Cref{thm:CF-recall} contained in \cite{torsors-regular}*{Proposition 5.3.6} or in \cite{Fed22a}*{Theorem 6} would suffice for us as well, but we prefer to take a clean general statement as our point of departure.

\bthm \label{thm:CF-recall}
Let $G$ be a reductive group over a ring $A$ and let $\sE$ be a $G$-torsor over $\bP^1_A$. 
\benum
\m \label{m:CFR-a}
{\upshape(\cite{totally-isotropic}*{Theorem 3.6}).} If $A$ is semilocal, then $\sE|_{\{t = 0\}} \simeq \sE|_{\{t = \infty\}}$.

\m \label{m:CFR-b}
{\upshape(\cite{totally-isotropic}*{Theorem 4.2}).} If $G$ is totally isotropic and $\sE|_{\{t = \infty\}}$ is trivial, then $\sE|_{\bA^1_A}$ is trivial.
\eenum
\ethm

\bpf
The claims are proved in a self-contained manner in the indicated references, although for \ref{m:CFR-a} we could alternatively cite \cite{PS24}. Let us briefly indicate what goes into the arguments. 

The key geometric input is the open immersion $i \colon {\mathbf B} G \hra \mathrm{Bun}_G$ from the algebraic $A$-stack ${\mathbf B} G$ parametrizing $G$-torsors over $A$ to the algebraic $A$-stack $\mathrm{Bun}_G$ parametrizing $G$-torsors over $\bP^1_A$, which one argues by using deformation theory for $G$-torsors. Moreover, in \ref{m:CFR-b} one uses Quillen patching for $G$-torsors over $\bA^1_A$ to reduce to local $A$. In both \ref{m:CFR-a} and \ref{m:CFR-b}, the geometry of $\mathrm{Bun}_G$ and the study of multiplicative group gerbes over $\bP^1_A$ allows one to pass to simply connected $G$. 

In both \ref{m:CFR-a} and \ref{m:CFR-b}, one knows the conclusion when $A$ is a field $k$ thanks to the classification of $G$-torsors over $\bP^1_k$ from, for instance, \cite{Ans18}, and the goal is to pass to semilocal $A$ using the open immersion $i$. This bootstrap is based on the Borel--Tits theorem \cite{Gil09}*{fait 4.3, lemme 4.5} (which uses the total isotropicity and the simply connectedness of $G$), by which certain glueings of trivial torsors can be obtained using ``elementary matrices.'' Since elementary matrices, and so the relevant glueings, lift across surjections, in \ref{m:CFR-b} one gets that  the $G$-torsor $\sE|_{\bA^1_A}$ extends to a $G$-torsor $\wt{\sE}$ over $\bP^1_A$ whose closed $A$-fiber is trivial; thanks to the openness of $i$, this means that $\wt{\sE}$, and so also $\sE|_{\bA^1_A}$, is trivial. The argument for \ref{m:CFR-a} is similar, except that, since $G$ is not totally isotropic, the lifting of glueings now happens along an $A$-(finite \'{e}tale) closed $Y \subset \bG_{m,\, A}$ such that $G_Y$ is totally isotropic. 
\epf



The following consequence of \Cref{thm:CF-recall}~\ref{m:CFR-b} is sharp in that it fails if the reductive $A$-group $G$ is no longer totally isotropic, see \cite{Fed16}*{Theorem~3 and what follows}.

 \bcor \label{conj:Horrocks}
 For a totally isotropic reductive group $G$ over a ring $A$ and an $A$-finite closed $Z \subset \bA^d_A$ with $d > 0$, every $G$-torsor over $\bA^d_A$ that trivializes over every affine $\bA^d_A \setminus Z$-scheme is trivial. 
\ecor



\bpf
Let $E$ be the $G$-torsor over $\bA^d_A$ in question. To show that $E$ is trivial, it suffices to show that its pullback under any section $s \in \bA^d_A(A)$ is trivial: indeed, as Gabber pointed out,  by applying this after base change to the coordinate ring $A[t]$ of $\bA^1_A$ and to the ``diagonal'' section of $\bA^1_{A[t]} \ra \Spec(A[t])$, we would get that $E$ itself is trivial.  Any $A$-point $s$ of $\bA^d_A$ factors through some $\bA^{d - 1}_A$-point, so we may replace $A$ by $A[t_1, \dotsc, t_{d - 1}]$ to reduce to $d = 1$. In the case $d = 1$, since the coordinate ring of $Z$ is a finite $A$-module, some monic polynomial in $A[t]$ vanishes on $Z$, so  we may replace $Z$ by this vanishing locus to arrange that $\bA^1_A \setminus Z$ be affine. The advantage of this is that then $E$ is even trivial over $\bA^1_A \setminus Z$. We then patch $E$ with the trivial torsor over $\bP^1_A \setminus Z$ to extend $E$ to a $G$-torsor over $\bP^1_A$ whose fiber at $\{ t = \infty\}$ is trivial. By \Cref{thm:CF-recall}~\ref{m:CFR-b}, then $E$ itself is trivial, as desired. 
\epf

\brem \label{rem:large-d-trivial}
In \Cref{conj:Horrocks}, if $d > 1$ and if the $G$-torsor in question trivializes over all of $\bA^d_A\setminus Z$ (not merely over every affine $\bA^d_A\setminus Z$-scheme), then the conclusion is an immediate consequence of \cite{EGAIV4}*{Proposition~19.9.8} and holds for any affine $A$-group $G$ (that need not be reductive).
\erem


 
\section{Overcoming the finite field obstructions}

A part of the reason of why we are able to progress beyond the cases of the Nisnevich conjecture established in \cite{Fed22b} is that in the critical \Cref{lem:pump-Y} below we find a way to bypass the finite field obstruction that hinders the geometric approach to the Nisnevich conjecture over finite fields. Even though in the Nisnevich case this obstruction is significantly more delicate, we still start with Panin's ``finite field tricks'' that have been used in every paper about the finite field or unramified mixed characteristic cases of the Grothendieck--Serre conjecture to overcome the corresponding obstacle in that context, see \cite{split-unramified}*{Lemma~6.1} or earlier works of Panin and of~Fedorov.


 \bd
For a ring $A$, a quasi-finite $A$-scheme $Z$, and an $A$-scheme $X$, there is no \emph{finite field obstruction} to embedding $Z$ into $X$ if for each maximal ideal $\fm \subset A$ with $k_\fm$ finite, we have
\be \label{eqn:ff-obstruction} \tag{$\dagger$}
\qq \#\{ z \in Z_{k_\fm} \,\vert\, [k_z : k_\fm] = m \} \le \#\{ z \in X_{k_\fm} \,\vert\, [k_z : k_\fm] = m  \} \qxq{for every} m \ge 1. 
\ee
\ed

\bprop \label{lem:pump-Y}
Let $A$ be a semilocal ring, let $Z$ be a quasi-finite, separated $A$-scheme, let $Y \subset Z$ be an $A$-finite closed subscheme, and let $X$ be an $A$-scheme such that for every maximal ideal $\fm \subset A$ with $k_\fm$ finite, some subscheme of $X_{k_\fm}$ is of finite type over $k_\fm$, positive dimensional, and geometrically irreducible. Suppose that $Y = Y_0 \sqcup Y_1$ with a $Y_0$ that has no finite field obstruction to embedding it into $X$. For every $n > 0$ and every large $N > 0$, there is a finite \'{e}tale surjection
\be\label{eqn:Ztilde-formula}
\wt{Z} \cong \underline{\Spec}(\sO_Z[t]/(f(t))) \surjects Z
\ee
with $f(t)$ monic of degree $N$ such that there is no finite field obstruction to embedding $\wt{Z}$ into $X$ and 
\[
\qqq \wt{Y} \ce Y \times_Z \wt{Z} \qxq{is a disjoint union}  \wt{Y} = \wt{Y}_0 \sqcup \wt{Y}_1 \qxq{such that}  \wt{Y}_0 \isomto Y_0
\]
and each connected component of $\wt{Y}_1$ is a scheme over $\Spec B$ for some finite $\bZ$-algebra $B$ each of whose residue fields $k$ of characteristic $p\mid n$ satisfies
\[
\#k > n\cdot \deg(\wt{Z}/Z).
\] 
\eprop

To be clear, the $\bZ$-algebra $B$ depends on the connected component of $\wt{Y}_1$ in question. 

\bpf
We may replace $Z$ by any $A$-finite scheme containing $Z$ as an open, so we use the Zariski Main Theorem \cite{EGAIV4}*{Corollaire 18.12.13} to assume that $Z = \Spec(A')$ for an $A$-finite $A'$.  To explain the role of the assumption on $X$, recall that by the Weil conjectures \cite{Poo17}*{Theorem 7.7.1~(ii)}, it implies that for every $d > 0$, every maximal ideal $\fm \subset A$ with $k_\fm$ finite, and every large $m > 0$,
\be \label{eqn:X-bound}
\#\{z \in X_{k_\fm}\, |\, [k_z : k_\fm] = m\} \ge d \qxq{(that is,}\tst \varinjlim_{m \rightarrow \infty} \#\{z \in X_{k_\fm}\, |\, [k_z : k_\fm] = m\} = \infty). 
\ee
Moreover, if the claim holds for $n$, then it also holds for every divisor of $n$ (with the same $\wt{Z}$). Thus, we may replace $n$ by any of its multiples, so we may assume that $n > 1$ and that it is divisible by all the positive residue characteristics of $A$. Moreover, we may assume that $Y$ contains all the closed points of $Z$ by adding some of these points to $Y_1$ if needed. Granted this, for each $N > 2$ we choose
\bitem
\m
an $f_{Y_0}(t) \in \bZ[t]$ that is the product of $t$ and a monic polynomial of degree $N - 1$ whose reduction modulo every prime $p \mid n$ is irreducible (and not linear because $N > 2$);

\m
a monic $f_{Y_1}(t) \in \bZ[t]$ of degree $N$ whose reduction modulo every prime $p \mid n$ is irreducible.

\eitem
We write $Y_i = \Spec(A'_i)$, view $f_{Y_i}(t)$  as an element of $A'_i[t]$, and choose  a monic polynomial $f(t) \in A'[t]$ whose image in $A'_i[t]$ is $f_{Y_i}(t)$. With $f(t)$ fixed, we let $\wt{Z}$ be defined by the formula \eqref{eqn:Ztilde-formula}. Since $f(t)$ is monic, this $\wt{Z}$ is finite and flat over $A$. To then check that $\wt{Z}$ is even finite \'etale over $A$ it suffices to check that the reduction of $f(t)$ modulo every maximal ideal $\fm \subset A$ is a separable polynomial over $k_\fm$. This is so by construction because $Y$ contains all the closed points of $Z$ and the images of the $f_{Y_i}(t)$ in $\mathbb{F}_p[t]$ with $p \mid n$ and also in $\mathbb{Q}[t]$ are separable (in fact, even either irreducible or a product of $t$ and a nonlinear irreducible polynomial). 

We let $\wt{Y}_0$ be the component of $Y_0 \times_Z \wt{Z}$ cut out by the factor $t$ of $f_{Y_0}(t)$, so that $\wt{Y}_0 \isomto Y_0$. By the choice of the $f_{Y_i}(t)$, each connected component of the complement $\wt{Y}_1$ of $\wt{Y}_0$ in $Y \times_Z \wt{Z}$ is an algebra over a finite $\bZ$-algebra $B$ that is either $\bZ[t]/(t^{-1}f_{Y_0}(t))$ or $\bZ[t]/(f_{Y_1}(t))$. Each residue field $k$ of characteristic $p > 0$ with $p \mid n$ of this $B$ has degree  either $N - 1$ or $N$ over $\bF_p$ and, for large $N$,
\[
\#k > nN = n\cdot \deg(\wt{Z}/Z).
\] 
It remains to show that there is no finite field obstruction to embedding $\wt{Z}$ into $X$. 
 An irreducible polynomial in $\bF_p[t]$ of degree $N$ splits into at most $i$ irreducible factors in $\bF_{p^i}[t]$, each of degree at least $N/i$. We now let $i$ range over the degrees of the finite residue fields of $Z$. By construction of $\wt{Z}$, we therefore get that, as $N$ grows, the number of closed points of $\wt{Z}$ not in $\wt{Y}_0$ with a finite residue field remains bounded by the sum of the degrees of the finite residue fields of $Z$. Moreover, as $N$ grows, the degrees of the finite residue fields of closed point of $\wt{Z}$ not in $\wt{Y}_0$ are all $\ge \eps N$ for some $\eps > 0$ that does not depend on $N$ (roughly, $\eps$ is the inverse of the maximum of the degrees of the finite residue fields of $Z$, except that we have to take it slightly smaller than that and let $N$ be large because the degree of $t^{-1}f_{Y_0}(t)$ is $N - 1$ and not $N$). In particular, for large $N$, by \eqref{eqn:X-bound}, there is no finite field obstruction to embedding the resulting $\wt{Z}$ into $X$: indeed, when $N$ is large, \eqref{eqn:ff-obstruction} with $\wt{Z}$ in place of $Z$ is automatic for $m < \epsilon N$ because there is no finite field obstruction to embedding $Y_0$ (so also $\wt{Y}_0$) into $X$ and $\wt{Y}_1$ does not contribute to the left side of \eqref{eqn:ff-obstruction}, whereas if $m \ge \epsilon N$, then the left side of \eqref{eqn:ff-obstruction} remains bounded while the right side tends to infinity in the view of \eqref{eqn:X-bound}.
\epf

\brems
\remi
As its proof shows, \Cref{lem:pump-Y} simplifies when $A$ is an $\bF_p$-algebra: then $B$ may be chosen to be a product of finite field extensions $k$ of $\bF_p$, each satisfying $\#k > n \cdot \deg(\wt{Z}/Z)$. 

\remi \label{rem:spread-out}
The $A$-quasi-finite $Z$ to be modified as in \Cref{lem:pump-Y} to avoid the finite field obstruction to embedding it into $X$ often occurs as a closed subscheme of  a smooth affine $A$-scheme $C$, and it is useful to lift the resulting  $\wt{Z} \surjects Z$ to a finite \'{e}tale cover $\wt{C} \surjects D$ of an affine open neighborhood $D \subset C$ of $Z$. Since $\wt{Z}$ is explicit, this is possible to arrange: it suffices to lift  $f(t)$ to a monic polynomial with coefficients in the coordinate ring of the semilocalization of $C$ at the closed points of $Z$ (built via prime avoidance \cite{SP}*{\href{https://stacks.math.columbia.edu/tag/00DS}{Lemma 00DS}}) and to spread out. 
\erems

The absence of finite field obstructions lets us reembed finite schemes $Z$ into $\bA^1_A$ as follows. This reembedding statement extends \cite{totally-isotropic}*{Lemma 2.5} and \cite{split-unramified}*{Lemma 6.3} (so also earlier versions due to Panin and Fedorov, see \emph{loc.~cit.}), but for applications to the Nisnevich conjecture we critically need its aspect about the compatibility $f|_Y = \iota_Y$, which the previous references do not supply. 


 
 

 \bprop \label{prop:embed-into-A1}
Let $A$ be a semilocal ring, let $U \subset \bA^1_A$ be an $A$-fiberwise nonempty open, and let $Z$ be a finite $A$-scheme. If there is no finite field obstruction to embedding $Z$ into $U$ and $Z$ is a closed subscheme of \emph{some} $A$-smooth affine scheme $C$ of relative dimension~$1$, then there is a closed immersion $\iota \colon Z \hra U$. Moreover, then $\iota$ may be chosen to be excisive\ucolon there are an affine open $D \subset C$ containing $Z$ and an \'{e}tale $A$-morphism $f\colon D \ra U$ that fits into a Cartesian square
\be \ba \label{eqn:EIA-square}
\xymatrix{
Z \ar@{^(->}[r] \ar[d]_-{\sim} & D \ar[d]_-{f}\\ 
Z' \ar@{^(->}[r] & U,
}
\ea\ee
in particular, such that $f$ embeds $Z$ as a closed subscheme $Z' \subset U$\uscolon in addition, for every $A$-finite closed subscheme $Y \subset Z$ and an embedding $\iota_Y \colon Y \hra U$, there are $D$ and $f$ as above with $f|_Y = \iota_Y$.
\eprop

\bpf
We fix embeddings $Z \subset C$ and $\iota_Y \colon Y \hra U$, let $\eps_Z \subset C$ be the first infinitesimal neighborhood of $Z$ in $C$, so that $\eps_Z$ is also finite over $A$, and let $k$ be the product of the residue fields of the maximal ideals of $A$. Since there is no finite field obstruction to embedding $(\eps_Z)_k$ into $U$, by \cite{totally-isotropic}*{Lemma~2.4}, there is such an embedding $\wt{\iota}_{k} \colon Z_k \hra U_k$ that extends $(\iota_Y)_k$ on $Y_k$. The closed immersions $\wt{\iota}_k$ and $\iota_Y$ are compatible, so there is a global section of $\eps_Z$ whose restriction to $(\eps_Z)_k$ (resp.,~$Y$) is the $\wt{\iota}_k$-pullback (resp.,~$\iota_Y$-pullback) of the standard coordinate of $\bA^1_A$. By sending the standard coordinate of $\bA^1_A$ to this global section, we obtain an $A$-morphism $\wt{\iota}\colon \eps_Z \hra U$ that extends the fixed $\iota_Y$. By construction and the Nakayama lemma \cite{SP}*{Lemma \href{https://stacks.math.columbia.edu/tag/00DV}{00DV}}, this $\wt{\iota}$ is a closed immersion. Its restriction to $Z$ is then the desired closed immersion $\iota \colon Z \hra U$. 

By lifting the $\wt{\iota}$-pullbacks of the standard coordinate of $\bA^1_A$, we may extend $\wt{\iota}$ to an $A$-morphism $\wt{f} \colon C \ra \bA^1_A$. By construction, the \emph{a priori} open locus of $C$ where $\wt{f}$ is quasi-finite (see \cite{SP}*{Lemma~\href{https://stacks.math.columbia.edu/tag/01TI}{01TI}}) contains the points of $Z$. Thus, since $Z$ has finitely many closed points, we may use prime avoidance \cite{SP}*{Lemma \href{https://stacks.math.columbia.edu/tag/00DS}{00DS}} to shrink $C$ around $Z$ to arrange that $\wt{f}$ is quasi-finite. The flatness criteria \cite{EGAIV2}*{Proposition~6.1.5} and \cite{EGAIV3}*{Corollaire 11.3.11} then ensure that $\wt{f}$ is flat at the points of $Z$, so, by construction, $\wt{f}$ is even \'{e}tale at the points of $Z$. Consequently, we may shrink $C$ further around $Z$ to arrange that $\wt{f}$ is \'{e}tale and factors through $U$. A section of a separated \'{e}tale morphism, such as $\wt{f}\i(\wt{f}(Z)) \ra \wt{f}(Z)$, is an inclusion of a clopen subset, so, by shrinking $C$ around $Z$ once more, we arrange that $Z = \wt{f}\i(\wt{f}(Z))$. This equality means that the square \eqref{eqn:EIA-square} is Cartesian, so, granted all the shrinking above, it remains to set $D \ce C$ and $\wt{f} \ce f$. 
 \epf

\section{Grothendieck--Serre for smooth relative curves over arbitrary rings} \label{sec:abstract-GS}

We use the reembedding techniques discussed in \Cref{prop:embed-into-A1} to present a Grothendieck--Serre phenomenon over arbitrary base rings: in \Cref{thm:abstract-GS} we show that torsors under reductive groups over smooth relative curves are Zariski semilocally trivial as soon as they are trivial away from some relatively finite closed subscheme. To approach this beyond constant $G$, we first establish \Cref{lem:equate-groups} about equating reductive groups, which is a variant of \cite{PSV15}*{Theorem~3.6} of Panin--Stavrova--Vavilov and combines ideas from \cite{split-unramified}*{Lemma 5.1} with those from the survey \cite{torsors-regular}*{Chapter~6}. 

\bd[\cite{torsors-regular}*{($\star$) in the beginning of Section 6.2}] \label{def:star}
For a ring $A$ and an ideal $I \subset A$, we consider the following property of a set-valued functor $\sF$ defined on the category of $A$-algebras: 
\be \tag{$\bigstar$} \label{eqn:condition-star}
\ba
&\x{for every $x \in \sF(A/I)$, there are a faithfully flat, finite, \'{e}tale $A$-algebra $\wt{A}$,}\\
&\x{an $A/I$-point $a \colon \wt{A} \surjects A/I$, and an $\wt{x} \in \sF(\wt{A})$ whose $a$-pullback is $x$.}
\ea
\ee 
\ed

Of course, since $\wt{A}$ is \'etale over $A$, faithful flatness amounts to the surjectivity of $\Spec(\wt{A}) \ra \Spec(A)$. Moreover, any $\sF$ that is representable by a faithfully flat, finite, \'etale $A$-scheme satisfies \eqref{eqn:condition-star}.

\brem \label{rem:star-extension}
Let $f\colon \sF \ra \sF'$ be a map of functors on the category of $A$-algebras and, for a $y \in \sF'(A)$, let $\sF_y \subset \sF$ denote the $f$-fiber of $y$. If $\sF'$ has property \eqref{eqn:condition-star} with respect to $I \subset A$ and, for every faithfully flat, finite, \'{e}tale $A$-algebra $\wt{A}$ and every $y \in  \sF'(\wt{A})$, the fiber $(\sF|_{\wt{A}})_y$ has property \eqref{eqn:condition-star} with respect to any ideal $\wt{I} \subset \wt{A}$ with $\wt{A}/\wt{I} \cong A/I$, then $\sF$ itself has property \eqref{eqn:condition-star} with respect to $I \subset A$. This straight-forward d\'{e}vissage is useful in practice for dealing with short exact sequences.
\erem

\blem \label{lem:equate-groups}
For a semilocal ring $A$, an ideal $I \subset A$, reductive $A$-groups $G$ and $G'$ that on geometric $A$-fibers have the same type and whose maximal central tori $\rad(G)$ and $\rad(G')$ are isotrivial, maximal $A$-tori $T \subset G$ and $T' \subset G'$, and an $A/I$-group isomorphism 
\[
\iota \colon G_{A/I} \isomto G'_{A/I} \qxq{such that} \iota(T_{A/I}) = T'_{A/I}, 
\]
there are a faithfully flat, finite, \'{e}tale $A$-algebra $\wt{A}$ equipped with an $A/I$-point $a \colon \wt{A} \surjects A/I$ and an $\wt{A}$-group isomorphism $\wt{\iota} \colon G_{\wt{A}} \isomto G'_{\wt{A}}$ whose $a$-pullback is $\iota$ and such that $\wt{\iota}(T_{\wt{A}}) = T'_{\wt{A}}$. 
\elem

\bpf
By passing to connected components, we may assume that $\Spec(A)$ is connected, so that the types of the geometric fibers of $G$ and $G'$ are constant. The claim is that the functor 
\[
X \ce \underline{\Isom}_\gp((G, T), (G', T'))
\]
that parametrizes those group scheme isomorphisms between base changes of $G$ and $G'$ that bring $T$ to $T'$ has property \eqref{eqn:condition-star} with respect to $I \subset A$. By \cite{SGA3IIInew}*{expos\'{e} XXIV, corollaires~1.10 et 2.2~(i)}, the normalizer $N_{G^\ad}(T^\ad)$ of the $A$-torus $T^\ad \subset G^\ad$ induced by $T$ acts freely on $X$ and, thanks to the assumption about the geometric fibers of $G$ and $G'$, the quotient 
\[
\ov{X} \ce X/N_{G^\ad}(T^\ad)
\]
is a faithfully flat $A$-scheme that becomes constant \'{e}tale locally on $A$. We claim that $\ov{X}$ has property \eqref{eqn:condition-star} with respect to $I \subset A$, more generally, that each quasi-compact subset of $\ov{X}$ is contained is some $A$-(finite \'{e}tale) clopen subscheme of $\ov{X}$ (such a clopen satisfies \eqref{eqn:condition-star}, as we pointed out after \Cref{def:star}). The advantage of the claim about the existence of an $A$-(finite \'etale) clopen is that it suffices to argue it after base change along any finite \'{e}tale cover of $A$. Thus, we may combine our assumption on $\rad(G)$ and $\rad(G')$ with \cite{SGA3IIInew}*{expos\'{e} XXIV, th\'{e}or\`{e}me 4.1.5} to assume that both $G$ and $G'$ are split. In this case, however, \cite{SGA3IIInew}*{expos\'{e} XXIV, th\'{e}or\`{e}me 1.3 (iii) et corollaire~2.2~(i)} ensure that $\ov{X}$ is a constant $A$-scheme, so the claim is clear.

With the property \eqref{eqn:condition-star} of $\ov{X}$ in hand, by \Cref{rem:star-extension}, we may replace $A$ by a finite \'{e}tale cover to reduce to showing that every $N_{G^\ad}(T^\ad)$-torsor has property \eqref{eqn:condition-star}. However, $N_{G^\ad}(T^\ad)$ is an extension of a finite \'{e}tale $A$-group scheme by $T^\ad$ (see, for instance, \cite{torsors-regular}*{Section 1.3.2}), so we may repeat the same reduction based on \Cref{rem:star-extension} and be left with showing that every $T^\ad$-torsor has property \eqref{eqn:condition-star} with respect to $I \subset A$. By \cite{SGA3IIInew}*{expos\'e XXIV, th\'eor\`eme 4.1.5 (i)}, the assumed isotriviality of $\rad(G)$ ensures that the maximal torus $T^\ad \subset G^\ad$ is isotrivial, and hence, $A$ being semilocal, that every $T^\ad$-torsor over $A$ is isotrivial as well. The desired property \eqref{eqn:condition-star} for $T^\ad$-torsors then follows from \cite{torsors-regular}*{Corollary~6.3.2 and its proof} (based on building an equivariant projective compactification of the $A$-torus $T^\ad$ using toric geometry): indeed, although the statement of \emph{loc.~cit.}~assumes that the local rings of $A$ are geometrically unibranch, its proof uses this assumption only to ensure that both $T^\ad$ and its torsor in question are isotrivial, which we have argued directly, whereas the Noetherianity assumption may be arranged by a limit argument.
\epf

\brem \label{rem:equate-Borels}
\Cref{lem:equate-groups} continues to hold if instead of the maximal $A$-tori $T \subset G$ and $T' \subset G'$, the groups $G$ and $G'$ come equipped with fixed quasi-pinnings extending Borel $A$-subgroups $B \subset G$ and $B' \subset G'$, and if $\iota$ and $\wt{\iota}$ are required to respect these quasi-pinnings, see \cite{split-unramified}*{Lemma~5.1}.
\erem

We are ready for the following promised Grothendieck--Serre type result over arbitrary base rings. 

\bthm \label{thm:abstract-GS}
Let $A$ be a ring, let $B$ be an $A$-algebra, let $C$ be a smooth affine $A$-scheme of pure relative dimension $d > 0$, let $\sG$ be a totally isotropic reductive $(C \tensor_A B)$-group scheme that descends to a reductive $C$-group $\wt{\sG}$ whose maximal central torus $\rad(\wt{\sG})$ is isotrivial Zariski semilocally on $C$ \up{resp.,~that descends to a reductive $B$-group  $G$}, and let $\sP \subset \sG$ be a parabolic $(C \tensor_A B)$-subgroup that descends to a parabolic $C$-subgroup $\wt{\sP} \subset \wt{\sG}$ \up{resp.,~to a parabolic $B$-subgroup $P \subset G$}. Every $\sG$-torsor $\sE$ over $C \tensor_A B$ whose restriction to $(C\setminus Z) \tensor_A B$ for some $A$-finite $Z \subset C$ reduces to a $\sR_u(\sP)$-torsor  trivializes Zariski semilocally on $C$, that is, for every $c_1, \dotsc, c_n \in C$, there is an affine open $C'\subset C$ containing all the $c_i$ such that $\sE$ trivializes over~$C'\tensor_A B$. 
\ethm

\bpf
Let $A'$ be the semilocal ring of $C$ at $c_1, \dotsc, c_n$, so that, by a limit argument, it suffices to show that $\sE$ trivializes over $A' \tensor_A B$. After base change to $A'$ the map $\Spec(A') \ra C$ induces a ``diagonal'' section of $C$, so, by performing such a base change and replacing $B$ by $A'\tensor_A B$, we reduce to showing that, when $A$ is semilocal, the pullback of $\sE$ under $s \tensor_A B$ for any $s \in C(A)$ is trivial. In addition, we enlarge $Z$ if necessary to ensure that $s \in Z(A)$. 

Granted this reformulation of the goal statement, let $k$ be the product of the residue fields of the maximal ideals of $A$. It follows from the presentation lemma \cite{split-unramified}*{Proposition 3.6~(vii)} (choose $Y = \emptyset$ there), alternatively, from \Cref{lem:GQ-mixed-char} below (choose $\cO = k$ and $Z = \emptyset$ there), that there are a principal affine open $C' \subset C$ containing $Z_k$ and a smooth $k$-morphism $\pi_k\colon C'_k \ra \bA^{d -1}_k$ of pure relative dimension $1$. By lifting the images of the standard coordinates, $\pi_k$ lifts to a morphism $\pi\colon C' \ra \bA^{d - 1}_A$. By the fibral criterion \cite{EGAIV3}*{th\'eor\`eme~11.3.10}, this $\pi$ is flat, so even smooth of pure relative dimension 1, at every point of $Z_k$. Thus, by shrinking $C'$ while keeping $Z_k \subset C'$, so also $Z \subset C'$, we may arrange $\pi$ to be smooth. At this point, we may replace $C$ by $C'$ and $A$ by $A[t_1, \dotsc, t_{d - 1}]$ (so $B$ by $A[t_1, \dotsc, t_{d - 1}] \tensor_A B$) to reduce the initial statement to the case when $d = 1$. We then repeat the reductions of the paragraph above to make $A$ semilocal again, with an $s \in C(A)$. 

Granted the above reduction to $d = 1$ and the reformulation of the goal statement, we will reduce to the case when $\sG$ descends to a reductive $B$-group $G$, which, being the pullback of $\sG$ along $s \tensor_A B$, is totally isotropic, and $\sP \subset \sG$ descends to a parabolic $B$-subgroup $P \subset G$. For this, it suffices to focus on the case when $\sG$ lifts to a reductive $C$-group $\wt{\sG}$ for which $\rad(\wt{\sG})$ is isotrivial Zariski semilocally on $C$ and $\sP \subset \sG$ lifts to a parabolic $C$-subgroup $\wt{\sP} \subset \wt{\sG}$, and to reduce this case to when $\wt{\sG}$ descends to a reductive $A$-group $\wt{G}$ and $\wt{\sP} \subset \wt{\sG}$ descends to a parabolic $A$-subgroup $\wt{P} \subset \wt{G}$. We begin by defining the candidate $\wt{P} \subset \wt{G}$ simply as the $s$-pullback of $\wt{\sP} \subset \wt{\sG}$. 

By shrinking $C$ around the closed points of $Z$, we may assume that $\rad(\wt{\sG})$ is isotrivial, that $\wt{\sG}$ has a maximal torus $\wt{\sT} \subset \wt{\sG}$ defined over $C$ (see \cite{SGA3II}*{expos\'{e} XIV, corollaire 3.20}), and, by passing to clopens if needed, that the type of the geometric $C$-fibers of $\wt{\sG}$ is constant.  We let $\wt{T} \subset \wt{G}$ be the $s$-pullback of $\wt{\sT}$. By \Cref{lem:equate-groups} and spreading out, there are an affine open $D \subset C$ containing $Z$ and a finite \'{e}tale cover $\wt{C} \surjects D$ for which $s$ lifts to some $\wt{s} \in \wt{C}(A)$ such that $\wt{\sG}|_{\wt{C}} \simeq \wt{G}|_{\wt{C}}$ compatibly with the fixed identification of pullbacks along $\wt{s}$. Thus, we may replace $C$ and $s$ by $\wt{C}$ and $\wt{s}$, respectively, and reduce to the case when $\wt{\sG}$ descends, that is, to when $\wt{\sG} = \wt{G}_C$.  To now likewise descend $\wt{\sP}$, we first pass to clopens to assume that the type of $\wt{\sP}$ as a parabolic subgroup of $\wt{G}_{C}$ is constant on $C$. Then $\wt{P}_C$ and $\wt{\sP}$ are parabolic subgroups of $\wt{G}_C$ of the same type, so, by \cite{SGA3IIInew}*{expos\'{e} XXVI, corollaire 5.5 (iv)} and a limit argument, they are conjugate over some affine open neighborhood of $Z$ in $C$. Since parabolic subgroups are self-normalizing \cite{SGA3IIInew}*{expos\'{e} XXVI, proposition~1.2}, the $s$-pullback of a conjugating section lies in $\wt{P}$, so we may adjust by this $s$-pullback to make the conjugating section pull back to the identity by $s$. Thus, by shrinking $C$ and adjusting the identification between $\wt{\sG}$ and $\wt{G}_C$ by an aforementioned conjugation, we achieve the promised reduction to the case when $\wt{\sP} \subset \wt{\sG}$ descends to $\wt{P} \subset \wt{G}$. 

With $\sP \subset \sG$ now being the base change of $P \subset G$, we wish to reduce to the case when $C = \bA^1_A$. 
For this, we begin with our closed immersion $Z \hra C $ and combine \Cref{lem:pump-Y} (with $Y = Y_0$ there being the schematic image of our section $s$), Remark \ref{rem:spread-out}, and \Cref{prop:embed-into-A1} to reduce to when there is an \'{e}tale morphism $C \ra \bA^1_A$ and a Cartesian square 
\[
\xymatrix{
Z \ar@{^(->}[r] \ar@{=}[d] & C \ar[d]\\ 
Z \ar@{^(->}[r] & \bA^1_A.
}
\]
By \cite{SP}*{Lemma \href{https://stacks.math.columbia.edu/tag/01PG}{01PG}} applied to the quasi-coherent ideal sheaf of $Z\subset \bP^1_A$, the $A$-finite $Z \subset \bA^1_A$ is the scheme-theoretic intersection of $A$-finite, finitely presented closed subschemes of $\bA^1_A$ containing it. By a limit argument, our \'{e}tale map $C \ra \bA^1_A$ becomes an isomorphism already when based changed to a small enough some such closed subscheme. Thus, we may enlarge our $A$-finite $Z$ to make it finitely presented over $A$ while retaining the Cartesian square above. The square remains Cartesian after base change to $B$, so we may  apply excision for $\sR_u(P)$-torsors \cite{split-unramified}*{Lemma~7.2~(b), Example~7.3} (with a limit argument that reduces to the Noetherian setting of \emph{loc.~cit.}; facilitating this limit argument was the only purpose of making $Z$ finitely presented) and then use patching supplied, for instance, by \cite{torsors-regular}*{Proposition 4.2.1}, to descend $\sE$ to a $G$-torsor over $\bA^1_B$ whose restriction to $(\bA^1_A \setminus Z)\tensor_A B$ reduces to an $\sR_u(P)$-torsor. Effectively, we have reduced to the promised case~$C = \bA^1_A$. 

Once $C = \bA^1_A$, we may use the avoidance lemma \cite{split-unramified}*{Lemma 3.1} to enlarge our $A$-finite $Z \subset \bA^1_A$ to be the vanishing locus of some hypersurface in $\bP^1_A$, to the effect that $\bA^1_A \setminus Z$ becomes affine. Then \cite{SGA3IIInew}*{expos\'{e} XXVI, corollaire 2.2} ensures that  $\sE$ trivializes over $(\bA^1_A \setminus Z) \tensor_A B$. It then suffices to apply \Cref{conj:Horrocks} to conclude that $\sE$ is trivial.
\epf

\brem \label{rem:AB}
In the case when $B = A$, \Cref{thm:abstract-GS} holds even without assuming that $\sG$ is totally isotropic. Indeed, we only used the total isotropicity assumption in the very last sentence of the proof, in order to apply \Cref{conj:Horrocks}, and without it we could instead change coordinates to make $s$ be the section $t = 0$, extend $\sE$ to a $G$-torsor over $\bP^1_A$ by patching it with the trivial torsor at infinity, and conclude the desired triviality of $s^*(\sE)$ by applying \Cref{thm:CF-recall}~\ref{m:CFR-a} instead. 
\erem

\section{The mixed characteristic cases of our main result on the Nisnevich conjecture}

We deduce the mixed characteristic cases of \Cref{thm:main-DVR} from the Grothendieck--Serre phenomenon of \Cref{thm:abstract-GS}. To arrive at its relative curve setting, we use the following presentation lemma.


\blem[\cite{split-unramified}*{Proposition 4.1}] \label{lem:GQ-mixed-char}
For a smooth, affine scheme $X$ of relative dimension $d > 0$ over a semilocal Dedekind ring $\cO$, points $x_1, \dotsc, x_m \in X$, and a closed subscheme $Z \subset X$ of codimension $\ge 2$, there are an affine open $X' \subset X$ containing $x_1, \dotsc, x_m$, an affine open $S \subset \bA^{d - 1}_\cO$, and a smooth morphism $f \colon X' \ra S$ of relative dimension $1$ such that $X' \cap Z$ is $S$-finite. \QED
\elem

\brem \label{rem:GQ-equichar}
In the case when $\cO$ is a field, the same statement holds under the weaker assumption that $Z$ is merely of codimension $\ge 1$ in $X$, see \cite{split-unramified}*{Remark 4.3} or \Cref{lem:Gabber-Quillen} below (whose proof does not use any other results from the present article).
\erem

\bpp[The abstract maximal torus] \label{pp:abstract}
To every reductive group $G$ over a scheme $S$ one associates an $S$-torus $T_G$, the \emph{abstract maximal torus} of $G$ defined by \'{e}tale descent on $S$ as follows. \'{E}tale locally on $S$, the group $G$ has a Borel $B \subset G$, and, letting $\sR_u(B) \subset B$ denote the unipotent radical,~one~sets
\[
T_G \ce B/\sR_u(B).
\]
Up to a canonical isomorphism, this $T_G$ does not depend on the choice of $B$, and so it descends to the original $S$: indeed, any two Borels are Zariski locally conjugate and, up to multiplying by a section of $B$, the conjugating section is unique \cite{SGA3IIInew}*{expos\'{e}~XXVI, proposition 1.2, corollaire 5.2}, so it suffices to note that the conjugation action of $B$ on $T_G$ is trivial because the latter is abelian. 
\epp

\bpp[Proof of \Cref{thm:main-DVR}~\ref{m:mDVR-3}] \label{pp:pf-mixed-char}
We have a semilocal ring $R$ that is flat and geometrically regular over a Dedekind subring $\cO$, an $r \in \cO$, a reductive $R[\f1r]$-group $G$ that either extends to a quasi-split reductive $R$-group or descends to a quasi-split reductive $\cO[\f1r]$-group, and a generically trivial $G$-torsor $E$ over $R[\f1r]$. We need to show that $E$ is trivial, and we will do this by applying \Cref{thm:abstract-GS}.

We use Popescu theorem \cite{SP}*{Theorem \href{https://stacks.math.columbia.edu/tag/07GC}{07GC}} and a limit argument to reduce to the case when $R$ is a semilocal ring of a smooth affine $\cO$-scheme $X$. By passing to connected components if needed, we may assume that $X$ is connected, of constant relative dimension $d$ over $\cO$. If $d = 0$, then $R$, and so also $R[\f1r]$, is a semilocal Dedekind ring, and $E$ is trivial by \cite{Guo22a}*{Theorem 1}; therefore, we lose no generality by assuming that $d > 0$. By shrinking $X$ if needed, we may assume that $G$ (resp.,~$E$) begins life over $X$ (resp.,~over $X[\f1r]$). In the case when our original $G$ extends to a quasi-split reductive $R$-group, we shrink $X$ further to make $G$ extend to a quasi-split reductive $X$-group $\wt{G}$ and we fix  a Borel $X$-subgroup $B \subset \wt{G}$. In the case when our original $G$ over $R[\f1r]$ descends to a quasi-split $\cO[\f1r]$-group, we shrink $X$ further to make sure that our new $G$ over $X[\f1r]$ still descends to a quasi-split reductive $\cO[\f1r]$-group, and we fix a Borel $\cO[\f1r]$-subgroup $B$ of this descended group. 

By applying the valuative criterion of properness to $E/B_{X[\f1r]}$, we may choose an open $U \subset X[\f1r]$ with complement of codimension $\ge 2$ such that $E_U$ reduces to a generically trivial $B$-torsor $\cE^B$ over $U$. By purity for torsors under tori \cite{CTS79}*{corollaire 6.9}, the $T_G$-torsor $\cE^B/\sR_u(B)$ over $U$ extends to a generically trivial $T_G$-torsor over $X[\f1r]$. To proceed, we use the following claim.
\epp

\bcl \label{cl:abstract-T}
The abstract maximal torus of $G$ has no nontrivial generically trivial torsors over $R[\f1r]$:
\[
\tst H^1(R[\f1r], T_G) \hra H^1(\Frac(R[\f1r]), T_G). 
\]
\ecl

\bpf
By our assumption on $G$ and the base change compatibility of the formation of the abstract maximal torus of a reductive group (see \S\ref{pp:abstract}), our $(T_G)_{R[\f1r]}$ is the base change of a torus $\cT$ defined over a ring $A$ that is either $R$ or $\cO[\f1r]$. By \cite{CTS87}*{Proposition 1.3}, this $\cT$ has a flasque resolution
\[
0 \ra \cF \ra \Res_{A'/A}(\bG_m) \ra \cT \ra 0,
\]
where $A'$ is a finite \'{e}tale $A$-algebra and $\cF$ is a flasque $A$-torus. For now, all we need to know about flasque tori is that, by the regularity of $R[\f1r]$ and \cite{CTS87}*{Proposition 1.4, Theorem 2.2 (ii)},
\[
\tst H^2(R[\f1r], \cF) \hra H^2(\Frac(R[\f1r]), \cF). 
\]
This reduces our desired claim to the vanishing $\Pic(R[\f1r] \tensor_A A') \cong 0$, which we argue as follows. In the case $A = R$, the ring $A'$ is again regular semilocal, so every line bundle on $A'[\f1r]$ extends to a line bundle on $A'$, and hence is trivial, to the effect that $\Pic(A'[\f1r]) = 0$, as desired. In the case $A = \cO[\f1r]$, by \cite{Ser79}*{Chapter I, Section 4, Proposition 8}, the normalization of $\cO$ in $A'$ is a finite $\cO$-algebra $\cO'$, in particular, $\cO'$ is again a Dedekind ring. Thus, $R \tensor_\cO \cO'$ is a finite $R$-algebra, and hence is semilocal, but is also flat and geometrically regular over $\cO'$, so it is regular by \cite{SP}*{Lemma~\href{https://stacks.math.columbia.edu/tag/033A}{033A}}. Since $R[\f1r] \tensor_A A'$ is a localization of $R \tensor_\cO \cO'$, it again follows that $\Pic(R[\f1r] \tensor_A A') \cong 0$, as desired. 
\epf

Thanks to \Cref{cl:abstract-T}, we may shrink $X$ around $\Spec(R)$ to trivialize the $T_G$-torsor $\cE^B/\sR_u(B)$, in particular, to make $E_U$ reduce to an $\sR_u(B)$-torsor. Since the complement $X[\f1r] \setminus U$ is of codimension $\ge 2$, its closure $Z$ in $X$ is also of codimension $\ge 2$. Thus, by \Cref{lem:GQ-mixed-char}, we may shrink $X$ around $\Spec(R)$ to arrange that there exists an affine open $S \subset \bA^{d - 1}_\cO$ and a smooth morphism $f \colon X \ra S$ of relative dimension $1$ such that $Z$ is $S$-finite. We can now apply \Cref{thm:abstract-GS} with $A\ce \Gamma(S, \sO_S)$ and $B\ce A[\f1r]$ (and \S\ref{conv} for the isotriviality condition) to conclude that $E$ is trivial over $R[\f1r]$; of course, here we are crucially using our assumption that the element $r$ comes from the base ring $\cO$. 
\QED


\section{The relative Grothendieck--Serre conjecture} \label{sec:relative-GS}



In equal characteristic, the approach to \Cref{thm:main-DVR} is based on 
the following relative version of the Grothendieck--Serre conjecture  that is a mild improvement to \cite{Fed22a}*{Theorem~1} (with an earlier more restrictive case due to Panin--Stavrova--Vavilov \cite{PSV15}*{Theorem~1.1}). Its case \ref{m:RGS-ii}, included here for completeness, reproves the equal characteristic case of the Grothendieck--Serre conjecture.

\addtocounter{footnote}{-1}
\renewcommand{\thefootnote}{\fnsymbol{footnote}}

\bthm \label{thm:rel-GS}
For a regular semilocal ring $R$ containing a field $k$,\footnote{\emph{Added after publication.} One should assume, in addition, that $R$ is geometrically regular over $k$ (this is automatic if $k$ is perfect), skip the first paragraph of the proof, and directly apply the Popescu theorem as in the second paragraph of the proof. The reduction to the case when $k = \bF$ given in the first paragraph of the proof is incorrect: the base change of $W \otimes_k R$ along $a$ is $W \otimes_k (k \otimes_\bF R)$ where the parenthetical factor is viewed as a $k$-algebra via its second factor $R$, and this base change does not agree with $W \otimes_\bF R$, contrary to what is claimed in the proof.} a reductive $R$-group $G$, and an affine $k$-scheme $W$, no nontrivial $G$-torsor over $W \tensor_k R$ trivializes over $W \tensor_k \Frac(R)$ if either
\benumr
\m \label{m:RGS-i}
$G$ is totally isotropic\uscolon or

\m \label{m:RGS-ii}
if $W  \tensor_k R$ is semilocal, for instance, if $W = \Spec(k)$.
\eenum
\ethm

\addtocounter{footnote}{0}
\renewcommand{\thefootnote}{\arabic{footnote}}

\bpf
Let $E$ be a $G$-torsor over $W \tensor_k R$ that trivializes over $W \tensor_k \Frac(R)$, let $\bF \subset k$ be the prime subfield, and consider the $k$-algebra $k \tensor_\bF R$.  The composition $R \xra{a} k \tensor_\bF R \xra{b} R$, in which the second map uses the $k$-algebra structure of $R$, is the identity. The base change of $E$ along $\id_W \tensor_k a$ is a $G$-torsor over $W \tensor_\bF R$ that trivializes over $W\tensor_\bF \Frac(R)$. Thus, it suffices to settle the claim with $k = \bF$ because, by then base changing further along $\id_W \tensor_k b$, we would get the desired triviality~of~$E$. 

Since $k$ is now perfect, Popescu theorem \cite{SP}*{Theorem \href{https://stacks.math.columbia.edu/tag/07GC}{07GC}} expresses $R$ as a filtered direct limit of smooth $k$-algebras. Thus, by passing to connected components of $\Spec(R)$ and doing a limit argument, we may assume that $R$ is a semilocal ring of a smooth, affine, irreducible $k$-scheme $X$ of dimension $d \ge 0$ and that $G$ and $E$ are defined over all of $X$. Since $E$ trivializes over $W \tensor_k \Frac(X)$, is also trivializes over $W \times_k (X \setminus Z)$ for some closed $Z \subsetneq X$. If $d = 0$, then $E$ is trivial, and if $d > 0$, then we may apply the presentation lemma of \Cref{rem:GQ-equichar} to shrink $X$ around $\Spec(R)$ so that there exist an affine open $S \subset \bA^{d - 1}_k$ and a smooth morphism $X \ra S$ of relative dimension $1$ that makes $Z$ finite over $S$. With such a fibration into curves in hand, however, the triviality of $E$ over $W \tensor_k R$ is a special case of \Cref{thm:abstract-GS} (with \S\ref{conv} for the isotriviality condition) and \Cref{rem:AB} applied with $A = \Gamma(S, \sO_S)$ and $B = \Gamma(W \times_k S, \sO_{W \times_k S})$ in \ref{m:RGS-i}, and with $A = B = \Gamma(W \times_k S, \sO_{W \times_k S})$ in \ref{m:RGS-ii}.
\epf

We will apply \Cref{thm:rel-GS} with $W \subset \bA^1_k$, in which case we may sharpen the assumptions as follows.


\blem[\cite{Gil02}*{Corollaire 3.10}] \label{lem:Gille-input}
For a reductive group $G$ over a field $K$ and an open $U \subset \bP^1_K$, each generically trivial $G$-torsor $E$ over $U$ reduces to a torsor under a maximal $K$-split subtorus of $G$\uscolon in particular, if $U \subset \bA^1_K$, then, since $U$ has no nontrivial line bundles, $E$ is a trivial $G$-torsor. \QED
\elem

\addtocounter{footnote}{-1}
\renewcommand{\thefootnote}{\fnsymbol{footnote}}

\bcor \label{cor:check-axiom}
For a regular semilocal ring $R$ containing a field $k$,\footnote{\emph{Added after publication.} As in \Cref{thm:rel-GS}, one should assume, in addition, that $R$ is geometrically regular over $k$ (this is automatic if $k$ is perfect).} a totally isotropic reductive $R$-group $G$, and a nonempty open $W \subset \bA^1_k$, every generically trivial $G$-torsor over $W \tensor_k R$ is trivial.
\ecor

\addtocounter{footnote}{0}
\renewcommand{\thefootnote}{\arabic{footnote}}

\bpf
Thanks to \Cref{lem:Gille-input}, \Cref{thm:rel-GS}~\ref{m:RGS-i} applies and gives the desired triviality.
\epf

\section{Extending $G$-torsors over a finite \'{e}tale subscheme of a relative curve} \label{sec:extend-Y}

A crucial preparation to the equicharacteristic case of the Nisnevich conjecture is  a result about extending $G$-torsors over a finite \'{e}tale closed subscheme of a smooth relative curve that we deduce in \Cref{prop:extend-0} from the reembedding techniques of \Cref{prop:embed-into-A1}. For wider applicability, we present this extension result axiomatically---it loosely amounts to a reduction of the Nisnevich conjecture to the Grothendieck--Serre conjecture. The equicharacteristic relative Grothendieck--Serre conjecture settled in \Cref{thm:rel-GS} supplies the required axioms in our main case of interest.

\bd \label{def:excisive}
For a ring $A$, a contravariant, set-valued functor $F$ on the category of $A$-schemes that are complements of $A$-quasi-finite closed subschemes in smooth affine $A$-schemes of pure relative dimension $1$ is \emph{excisive} if for all Cartesian~squares 
\[
\q \xymatrix{
Z \ar@{^(->}[r] \ar[d]_-{\sim} & S \ar[d]^-f\\ 
Z' \ar@{^(->}[r] & S'
}
\]
in which the horizontal maps are closed immersions, $Z$ and $Z'$ are $A$-quasi-finite and finitely presented, $S$ and $S'$ are complements of $A$-quasi-finite closed subschemes in smooth affine $A$-schemes of pure relative dimension $1$, and $f$ is \'{e}tale and induces an indicated isomorphism $Z \isomto Z'$, we have 
\[
F(S') \surjects F(S) \times_{F(S \setminus Z)} F(S'\setminus Z').
\]
\ed

For instance, for a quasi-affine, flat, finitely presented $A$-group $G$, the functor $H^1(-, G)$  is excisive, see \cite{torsors-regular}*{Proposition 4.2.1}. The following lemma is critical for our argument for \Cref{thm:main-DVR}~\ref{m:mDVR-1}.

\blem \label{lem:change-to-Y}
Let $A$ be a ring, let $S$ be an $A$-scheme, let $Y \subset S$ be an $A$-\up{separated \'{e}tale} closed subscheme that is locally cut out by a finitely generated ideal, and consider the decomposition 
\[
Y \times_A Y = \Delta \sqcup Y'
\]
in which $\Delta \subset Y \times_A Y$ is the diagonal copy of $Y$. The following square is Cartesian\ucolon
\[
\xymatrix{
\Delta \ar@{^(->}[r] \ar[d]_-{\sim} & S_Y \setminus Y' \ar[d]\\ 
Y \ar@{^(->}[r] & S,
}
\]
in particular, if $F$ is an excisive functor as in Definition \uref{def:excisive} and $S$ is the complement of an $A$-quasi-finite closed subscheme in some smooth affine $A$-scheme of pure relative dimension $1$, then an element of $F(S\setminus Y)$ extends to $F(S)$ if and only if its pullback to $F((S \setminus Y)_Y)$ extends to $F(S_Y \setminus Y')$\uscolon for instance, for a quasi-affine, flat, finitely presented $S$-group $G$, a $G$-torsor over $S\setminus Y$ extends to a $G$-torsor over $S$ if and only if its base change to $(S \setminus Y)_Y$ extends to a $G$-torsor over $S_Y \setminus Y'$. 
\elem

\bpf
The claimed decomposition $Y \times_A Y = \Delta \sqcup Y'$ exists because any section of a separated \'{e}tale morphism, such as the projection $Y \times_A Y \ra Y$, is both a closed immersion and an open immersion. Thus, the square in question is Cartesian because the \'{e}tale map $S_Y \setminus Y' \ra S$ induces an isomorphism $\Delta \isomto Y$. The claim about $F$ is then immediate from \Cref{def:excisive}.
\epf

We are ready for our key axiomatic extension result, which extends Fedorov's \cite{Fed22b}*{Proposition~2.8}.

\bprop \label{prop:extend-0}
Let 
\bitem
\m
$A$ be a reduced semilocal ring that contains a field \up{so also a field $k$ that is either $\bQ$ or $\bF_p$}, 

\m
$C$ be a smooth affine $A$-scheme of pure relative dimension $1$, 



\m
$Y \subset C$ be an $A$-\up{finite \'{e}tale} closed subscheme, and

\m
$F$ be an excisive, pointed set valued functor as in Definition~\uref{def:excisive}.
\eitem
Suppose that for  each finite \'etale $k$-algebra $k^\prime$, each finite \'{e}tale $Y$-scheme $\cY$ that is also a $k^\prime$-scheme, each $\sY \subset \bA^1_\cY$ that is both a union of finitely many pairwise disjoint $\cY$-points and a base change of a finite set of $k^\prime$-points of $\bA^1_{k^\prime}$, and each $\cY$-finite closed subscheme $\sZ \subset \bA^1_{\cY}$ containing $\sY$, we have
\be \label{eqn:assumption-on-F-0}
\Ker(F(\bA^1_{\cY}) \ra F(\bA^1_{\cY}\setminus \sZ)) \surjects \Ker(F(\bA^1_\cY \setminus \sY) \ra F( \bA^1_\cY \setminus \sZ)),
\ee
that is, every element of $F(\bA^1_\cY \setminus \sY)$ that trivializes away from some $\cY$-finite $\sZ \subset \bA^1_\cY$ containing $\sY$ extends to $F(\bA^1_\cY)$. Then, for every $A$-finite closed subscheme $Z \subset C$ containing $Y$, we have
\be \label{eqn:conclusion-about-F-0}
\Ker(F(C) \ra F(C\setminus Z)) \surjects \Ker(F(C \setminus Y) \ra F(C \setminus  Z)),
\ee
that is, every element of $F(C \setminus Y)$ that trivializes away from some $A$-finite $Z \subset C$ containing $Y$ extends to $F(C)$.
\eprop

 \Cref{cor:check-axiom} supplies the assumption \eqref{eqn:assumption-on-F-0} when $A$ is regular of equicharacteristic and $F(-)$ is $H^1(-, G)$ for a reductive $A$-group $G$ such that $G_Y$ is totally isotropic.

For proving \Cref{prop:extend-0} and, simultaneously, for potential future applications in mixed characteristic, it is convenient to directly argue the following more general statement in \Cref{prop:extend}. It incorporates an auxiliary larger curve $C'$ to help with intermediate reductions in the proof and it also  works over $\mathbb{Z}$ instead of over a base field $k$. Since the finite \'etale $Y$-scheme $\cY$ in \Cref{prop:extend-0} is reduced, any map from a finite $\mathbb{Z}$-algebra $B$ to the coordinate ring of  $\cY$ factors through some $k^\prime$ as in \Cref{prop:extend-0}, so \Cref{prop:extend-0} is indeed a special case of \Cref{prop:extend}. 



\bprop \label{prop:extend}
Let 
\bitem
\m
$A$ be a semilocal ring, 

\m
$C'$ be a smooth affine $A$-scheme of pure relative dimension $1$, 

\m
$Y^\prime \subset C'$ is an $A$-\up{finite \'etale} closed subscheme with complement $C := C' \setminus Y'$,


\m
$Y \subset C$ be an $A$-\up{finite \'{e}tale} closed subscheme, and

\m
$F$ be an excisive, pointed set valued functor as in Definition~\uref{def:excisive}.
\eitem
Suppose that for each finite $\bZ$-algebra $B$, each finite \'{e}tale $(Y\cup Y')$-scheme $\cY$ that is also a $B$-scheme, each $\sY \subset \bA^1_\cY$ that is both a union of finitely many pairwise disjoint $\cY$-points and a base change of a finite set of \up{possibly nondisjoint} $B$-points of $\bA^1_B$, and each $\cY$-finite closed $\sZ \subset \bA^1_{\cY}$ containing $\sY$,
\be \label{eqn:assumption-on-F}
\Ker(F(\bA^1_{\cY}) \ra F(\bA^1_{\cY}\setminus \sZ)) \surjects \Ker(F(\bA^1_\cY \setminus \sY) \ra F( \bA^1_\cY \setminus \sZ)).
\ee
Then, for every $A$-finite closed subscheme $Z \subset C'$ containing~$Y \cup Y'$, we have
\be \label{eqn:conclusion-about-F}
\Ker(F(C) \ra F(C\setminus Z)) \surjects \Ker(F(C \setminus Y) \ra F(C \setminus  Z)).
\ee
\eprop

 The proof is a formal reduction of the property \eqref{eqn:conclusion-about-F} to the case when $C = \bA^1_A$ and $Y$ is ``constant.''

\bpf
We fix an $\gA \in \Ker(F(C\setminus Y) \ra F(C \setminus Z))$ that we wish to extend over $Y$. The assumption that $F$ be excisive is stable under finite \'etale base change in $A$. Thus, we may use \Cref{lem:change-to-Y} to base change along $Y \ra \Spec(A)$ and shrink the base changed $C$ by removing the off-diagonal part of $Y \times_A Y$ to reduce to the case when $Y \cong \Spec(A)$ (so the base changed $C'$ is kept and the base changed $Y'$ is enlarged by uniting it with the off-diagonal part of $Y \times_A Y$). Moreover, we decompose $A$ to reduce to the case when $\Spec(A)$ is connected, so that $\deg(Y \cup Y' /A)$ is a well-defined integer.   

We let $n$ be the product of $\deg((Y \cup Y')/A)$ and of all the prime numbers $p$ with $p \not\in A^\times$. By \Cref{lem:pump-Y} (applied with $Y_0 = Y$ and $Y_1 = Y'$) with Remark~\ref{rem:spread-out}, there are an affine open $D \subset C'$ containing $Z$ (so also $Y \cup Y'$) and a finite \'{e}tale cover $\wt{C}' \surjects D$ such that there is no finite field obstruction to embedding $\wt{Z} \ce Z \times_{C'} \wt{C}'$  into $\bA^1_A$ and
\[
\wt{Y} \ce Y \times_{C'} \wt{C}' \qxq{decomposes as}  \wt{Y} = \wt{Y}_0 \sqcup \wt{Y}_1 \qxq{such that} \wt{Y}_0 \isomto \Spec(A)
\]
and each component of $\wt{Y}_1$ or of $\wt{Y}' \ce Y' \times_{C'} \wt{C}' $  is a scheme over some finite $\bZ$-algebra $B$ each of whose residue fields $k$ of characteristic $p \mid n$ satisfies
\[
\#k > n \cdot \deg(\wt{Z}/Z) \ge \deg((\wt{Y} \cup \wt{Y}')/A).
\] 
By construction, setting $\wt{C} \ce (D \setminus Y') \times_{C'} \wt{C}'$, we have a Cartesian square
\[
\xymatrix{
\wt{Y}_0 \ar[d]_-{\sim} \ar@{^(->}[r] & \wt{C} \setminus \wt{Y}_1 \ar[d]  \\
Y \ar@{^(->}[r] & C \cap D. 
}
\]
Thus, since $F$ is excisive, to extend $\gA$ over $Y$ we may first restrict to $C\cap D$ (in \Cref{def:excisive}, choose $f$ to be the inclusion $C \cap D \hookrightarrow C$ and choose $Z$ and $Z'$ there to be our $Y$) and then pass to $\wt{C} \setminus \wt{Y}_1$. In other words, we may replace $Y \subset C \subset C'$ by $\wt{Y}_0 \subset \wt{C}\setminus \wt{Y}_1 \subset \wt{C}'$ and $\gA$ by its pullback to $\wt{C} \setminus \wt{Y}$ to reduce to the case when $Y \cong \Spec(A)$ and each connected component of $Y'$ is an algebra over  some finite $\bZ$-algebra $B$ each of whose residue fields $k$ of characteristic $p \mid n$ satisfies $\#k > \deg((Y \cup Y' )/A)$ and there is no finite field obstruction to embedding $Z$ into $\bA^1_A$  (the goal of this step is to prepare for reducing to $\mathbb{A}^1_A$ by excision afterwards).   By \Cref{prop:embed-into-A1}, such an embedding then exists, more precisely, there are an affine open $D \subset C'$ containing $Z$ and a Cartesian square
\be\ba \label{eqn:aux-Cartesian}
\xymatrix{
Z \ar@{^(->}[r] \ar[d]_-{\sim} & D \ar[d]_-{f}\\ 
Z' \ar@{^(->}[r] & \bA^1_A
}
\ea \ee
in which the map $f$ is \'{e}tale and embeds $Z$ as a closed subscheme $Z' \subset \bA^1_A$. The square remains Cartesian after passing to the complements of the $A$-(finite \'{e}tale) $Y \cup Y' $ viewed inside $Z$ (so also inside $Z'$). Thus, for the purpose of extending $\gA$ over $Y$, we may use the excisive property of $F$ to patch the restriction $\gA|_{D \setminus (Y \cup Y')}$ with the origin in $F(\bA^1_A\setminus Z')$ to reduce to the case when $C' = \bA^1_A$. 

In conclusion, at the cost of stepping back to the setting of a more general $Y$, we have reduced our overall sought claim about extending $\alpha$ to the case when $C' = C = \bA^1_A$ and $Y \cong \Spec(A) \sqcup y$ are such that each connected component of $y$ is a scheme over some finite $\bZ$-algebra $B$ each of whose residue fields $k$ of characteristic $p \mid n$ satisfies $\# k > \deg(Y/A)$. To extend $\gA$ over any fixed connected component of $y$, since $F$ is excisive, we may base change to this component and use \Cref{lem:change-to-Y} (noting that we may use the same $n$ after such a base change and that $\deg(Y/A)$ is stable under such a base change). Thus, we may assume that $A$ itself is an algebra over some finite $\bZ$-algebra $B$ as above: indeed, once we argue the claim under this assumption, by the previous sentence, we will be able to extend $\gA$ over $y$ by iteratively extending over each of its components, and this will 
leave us with the case $y = \emptyset$, in which case we may choose $B = \bZ$ to force the same assumption.

Granted the reductions above, we now induct on the number of disjoint copies of $\Spec(A)$ contained in $Y$ to reduce to when $Y \simeq \bigsqcup \Spec(A)$. Indeed, suppose that $Y$ has a connected component $W$ that does not map isomorphically to $\Spec(A)$, so that $W$ is of degree $\ge 2$ over $A$. Since $W \times_A W$ contains the diagonal copy of $W$ as a clopen (compare with \Cref{lem:change-to-Y}), the $W$-(finite \'{e}tale) closed subscheme $Y \times_A W \subset \bA^1_{W}$ has the same degree $\deg(Y/A)$ over $W$ and contains strictly more disjoint copies of $W$ than $Y$ contained disjoint copies of $\Spec(A)$. Thus, by the inductive hypothesis, the pullback of $\gA$ to  $\bA^1_{W} \setminus (Y \times_A W)$ extends over $Y \times_A W$. By \Cref{lem:change-to-Y}, this implies that $\gA$ extends over $W$. By repeating this for each possible $W$, we effectively eliminate connected components of $Y$ one by one until we reduce to the desired base case when $Y \simeq \bigsqcup \Spec(A)$. 

To treat this last case, we set $m \ce \deg(Y/A)$, so that, without losing generality, $m \ge 1$, and we will use our assumption \eqref{eqn:assumption-on-F}. We take $\cY \ce \mathrm{Spec}(A)$, which is, by our assumption, a $B$-scheme. 
However, we cannot simply choose $\sY = Y$ because, even though $Y$ is a union of $m$ pairwise disjoint $A$-points of $\mathbb{A}^1_A$, these points need not be defined over $B$, that is, $Y$ need not be a base change of a finite set of $B$-points of $\mathbb{A}^1_B$ (not even up to an automorphism of $\mathbb{A}^1_A$ if $m \ge 3$). Nevertheless, the condition on the residue fields of $B$ does ensure that $\mathbb{A}^1_B$ has $m$ distinct $B$-points that pull back to $m$ pairwise distinct $k$-points of
 $\mathbb{A}^1_k$ for every residue field $k$ of $B$ of characteristic $p$ with $p \not \in A^\times$. The union of these $B$-points of $\bA^1_B$ base changes to a closed subscheme $\sY \subset \bA^1_A$ that is a union of $m$ pairwise disjoint $A$-points of $\mathbb{A}^1_A$ (disjointness may be tested over the residue fields of the maximal ideals of $A$). This last condition ensures that there is an $A$-isomorphism $Y \simeq \sY$, and \Cref{prop:embed-into-A1} (especially, its final aspect) then supplies an affine open $D \subset \bA^1_A$ containing $Z$ and a Cartesian square as in \eqref{eqn:aux-Cartesian} such that $f$ maps $Y$ isomorphically onto $\sY$. Thus, since $F$ is excisive, we reduce to the case when $Y = \sY$ inside $\bA^1_A$. At this point we conclude by applying our assumption \eqref{eqn:assumption-on-F}. 
\epf

\addtocounter{footnote}{-1}
\renewcommand{\thefootnote}{\fnsymbol{footnote}}

\bcor \label{thm:A1}
Let $R$ be a regular semilocal ring containing a field, let $G$ be a totally isotropic reductive $R$-group scheme, let $C$ be a smooth affine $R$-scheme of pure relative dimension $1$, and let $Y \subset C$ be an $R$-\up{finite \'{e}tale} closed subscheme. Every $G$-torsor over $C \setminus Y$ that trivializes away from some $R$-finite closed subscheme $Z \subset C$ containing $Y$ extends to a $G$-torsor over $C$. 
\ecor

\bpf
By \Cref{cor:check-axiom}, for a product of fields $k^\prime$,\footnote{\emph{Added after publication.} One should restrict to $k^\prime$ that are products of perfect fields.} a $k^\prime$-fiberwise nonempty open $W \subset \mathbb{A}^1_{k^\prime}$, and a finite \'{e}tale $R$-algebra $R^\prime$ that is a $k^\prime$-algebra, every generically trivial $G$-torsor over $W \tensor_{k^\prime} R^\prime$ is trivial. Thus, the excisive functor $F(-) := H^1(-, G)$ fulfils the axiomatic assumption \eqref{eqn:assumption-on-F-0} (let $W \subset \mathbb{A}^1_{k^\prime}$ be such that $W_\cY = \mathbb{A}^1_\cY \setminus \sY$). In effect, \Cref{prop:extend-0} applies  and gives the claim.
\epf

\addtocounter{footnote}{0}
\renewcommand{\thefootnote}{\arabic{footnote}}

\section{The Nisnevich conjecture over a field}


The final preparation to the equicharacteristic case of the Nisnevich conjecture is the following geometric presentation lemma in the spirit of Gabber's refinement \cite{Gab94b}*{Lemma~3.1} of the Quillen presentation lemma \cite{Qui73}*{Section 7, Lemma 5.12}, which itself is a variant of the Noether normalization theorem. For us, it is crucial to have its aspect about the smooth divisor $D$.


\blem \label{lem:Gabber-Quillen}
For a smooth, affine, irreducible scheme $X$ of dimension $d > 0$ over a field $k$ that is either finite or of characteristic $0$,\footnote{The assumption on $k$ is likely not optimal but it will suffice and we do not wish to further complicate the proof.} points $x_1, \dotsc, x_m \in X$, 
a proper closed subscheme $Z \subset X$, and a $k$-smooth divisor $D \subset X$, there are an affine open $X' \subset X$ containing $x_1, \dotsc, x_m$, an affine open $S \subset \bA^{d - 1}_k$, and a smooth morphism
\[
f\colon X' \ra S
\]
of relative dimension $1$ such that 
\[
X' \cap Z = f\i(S) \cap Z \qxq{is $S$-finite and} X' \cap D = f\i(S) \cap D \qx{is $S$-\up{finite \'{e}tale}.}
\]
\elem

\bpf
In the case $d = 1$, we may choose $X' = X$ and $S = \Spec(k)$, so we assume that $d > 1$. We also replace each $x_i$ by a specialization to reduce to $x_i$ being a closed point (see \cite{SP}*{Lemma~\href{https://stacks.math.columbia.edu/tag/02J6}{02J6}}), and in this case we will force each $f(x_i)$ to be the origin of $\bA^{d - 1}_k$. We embed $X$ into some projective space $\bP^N_k$ 
and then form closures to arrange that $X$ is an open of a projective $\ov{X} \subset \bP^N_k$ of dimension $d$ with $\ov{X} \setminus X$ of dimension $\le d - 1$ and that there are
\bitem
\m
a projective $\ov{D} \subset \ov{X}$ of dimension $d - 1$ with $\ov{D} \setminus D$ of dimension $\le d - 2$, and

\m
a projective $\ov{Z} \subset \ov{X}$ of dimension $\le d - 1$ with $\ov{Z} \setminus Z$ of dimension $\le d - 2$.
\eitem
We use the avoidance lemma \cite{GLL15}*{Theorem 5.1} and postcompose with a Veronese embedding to build a hyperplane $H_0$ not containing any $x_i$ such that $(\ov{X} \setminus X) \cap H_0$ is of dimension $\le d - 2$ (to force the dimension drop, choose appropriate auxiliary closed points and require $H_0$ to not contain them). By the Bertini theorem \cite{Poo04}*{Theorem 1.3} of Poonen if $k$ is finite and by the Bertini theorem of \cite{split-unramified}*{second paragraph of the proof of Lemma 3.2} applied both to $X$ and to $D$ in place of $X$ if $k$ is of characteristic $0$, there is a hypersurface $H_1 \subset \bP^N_k$ such that
\bitem
\m
$H_1$ contains $x_1, \dotsc, x_m$;

\m
$X \cap H_1$ (resp.,~$D \cap H_1$) is $k$-smooth of dimension $d - 1$ (resp.,~$d - 2$);

\m
$Z\cap H_1$ is (resp.,~$(\ov{D} \setminus D)\cap H_1$ and $(\ov{Z} \setminus Z)\cap H_1$ are) of dimension $\le d - 2$ (resp.,~$\le d - 3$);

\m
$(\ov{X} \setminus X) \cap H_0 \cap H_1$ is of dimension $\le d - 2$.
\eitem
In particular, by passing to intersections with $H_1$, we are left with an analogous situation with $d$ replaced by $d - 1$. Therefore, by iteratively applying the Bertini theorem in this way, we build hypersurfaces $H_1, \dotsc, H_{d - 1}$ such that 
\benumr
\m \label{m:IC-o}
the $x_1, \dotsc, x_m$ lie in $H_1 \cap \dotsc \cap H_{d - 1}$ but not in $H_0$;

\m \label{m:IC-i}
$X \cap H_1 \cap \dotsc \cap H_{d - 1}$ (resp.,~$D \cap H_1 \cap \dotsc \cap H_{d - 1}$) is $k$-smooth of dimension $1$ (resp.,~$k$-\'{e}tale);

\m \label{m:IC-ii}
$(\ov{D} \setminus D) \cap H_1 \cap \dotsc \cap H_{d - 1} = (\ov{Z} \setminus Z) \cap H_1 \cap \dotsc \cap H_{d - 1} = \emptyset$.

\m \label{m:IC-iii}
$(\ov{X} \setminus X) \cap H_0 \cap H_1 \cap \dotsc \cap H_{d - 1} = \emptyset$.
\eenum
By letting $1, w_1, \dotsc, w_{d - 1}$ be the degrees of the hypersurfaces $H_0, H_1, \dotsc, H_{d - 1}$ and choosing defining equations $h_i$ of the $H_i$, we determine a projective morphism $\wt{f}\colon \wt{X} \ra \bP_k(1, w_1, \dotsc, w_{d - 1})$ from the weighted blowup $\wt{X} \ce \Bl(h_0, \dotsc, h_{d - 1})$ to the weighted projective space such that the diagram
\[
\xymatrix{
\ov{X} \setminus H_0 \ar@{^(->}[r] \ar[d]_{f} &\ov{X} \setminus (H_0 \cap \dotsc \cap H_{d - 1}) \ar@{^(->}[r] \ar[d] &\wt{X} \ar[d]_-{\wt{f}} \\ 
\bA_k^{d - 1}\ar@{^(->}[r]&\bP_k(1, w_1, \dotsc, w_{d - 1}) \ar@{=}[r]& \bP_k(1, w_1, \dotsc, w_{d - 1})
}
\]
commutes, where the bottom left arrow is the inclusion of the open locus where the first standard coordinate of $\bP_k(1, w_1, \dotsc, w_{d - 1})$ does not vanish, see \cite{split-unramified}*{Sections 3.4 and 3.5}. By \ref{m:IC-o}, each $f(x_i)$ is the origin of $\bA^{d - 1}_k$. By \ref{m:IC-i} and the dimensional flatness criterion \cite{EGAIV2}*{Proposition~6.1.5}, at every point of the fiber above the origin of $\bA^{d - 1}_k$, the map $f$ is smooth of relative dimension $1$ and its restriction to $D$ is \'{e}tale. Since $\wt{f}$ is projective, \ref{m:IC-ii}--\ref{m:IC-iii} and the openness of the quasi-finite locus \cite{SP}*{Lemma \href{https://stacks.math.columbia.edu/tag/01TI}{01TI}} ensure that for some affine open neighborhood of the origin $S \subset \bA^{d - 1}_k$ both $f\i(S) \cap Z$ and $f\i(S) \cap D$ are $S$-finite (see also \cite{SP}*{Lemma \href{https://stacks.math.columbia.edu/tag/02OG}{02OG}}). In conclusion, any affine open of $f\i(S)$ that contains all the $x_i$ and all the points of $Z$ and $D$ that lie above the origin of $\bA^{d - 1}_k$ becomes a sought $X'$ after possibly shrinking $S$ further.
\epf

\addtocounter{footnote}{-2}
\renewcommand{\thefootnote}{\fnsymbol{footnote}}

\bpp[Proof of \Cref{thm:main-DVR}~\ref{m:mDVR-1}] \label{pp:Nisnevich-conj-pf}
We have a regular semilocal ring $R$ containing a field $k$,\footnote{\emph{Added after publication.} One should immediately replace $k$ by its prime subfield to assume that $k$ is either $\bQ$ or $\bF_p$ (as at the end of this paragraph), and hence reduce to the case when $k$ is perfect.} a regular parameter $r \in R$, a reductive $R$-group $\cG$ with $\cG_{R/(r)}$ totally isotropic, and a generically trivial $\cG$-torsor $E$ over $R[\f{1}{r}]$. We need to show that $E$ is trivial, equivalently, by a known case of the Grothendieck--Serre conjecture \Cref{thm:rel-GS}~\ref{m:RGS-ii}, we need to extend $E$ to a $\cG$-torsor $\cE$ over $R$. For this, by Zariski patching and a limit argument, we may semilocalize $R$ along the union of those maximal ideals $\fm \subset R$ that contain $r$ and reduce ourselves to the case when $r$ lies in every maximal ideal $\fm \subset R$. Moreover, we may replace $k$ by its prime subfield to assume that $k$ is either~$\bQ$~or~some~$\bF_p$. 

Popescu theorem \cite{SP}*{Theorem \href{https://stacks.math.columbia.edu/tag/07GC}{07GC}} expresses $R$ as a filtered direct limit of smooth $k$-algebras. Thus, by passing to connected components of $\Spec(R)$ and doing a limit argument, we may assume that $R$ is a semilocal ring of a smooth, affine, irreducible $k$-scheme $X$ of dimension $d \ge 0$, that $r$ is a global section of $X$ that cuts out a $k$-smooth divisor $D \subset X$ with complement $U \ce X \setminus D$, that $\cG$ (resp.,~$E$) is defined over all of $X$ (resp.,~$U$), and that $\cG_D$ is totally isotropic. Since $E$ is trivial over $\Frac(X)$, there is a closed $\sZ \subsetneq X$ containing $D$ such that $E$ is trivial over $U \setminus \sZ$. If $d = 0$, then $E$ is trivial, so we assume that $X$ is of dimension $d > 0$. Finally, we use \cite{SGA3II}*{expos\'{e}~XIV, corollaire~3.20} to shrink $X$ further to make $\cG$ have a maximal torus $T$ defined over all of $X$.

With these preparations, \Cref{lem:Gabber-Quillen} allows us to shrink $X$ around $\Spec(R)$ to arrange that there exist an affine open $S \subset \bA^{d - 1}_k$ and a smooth morphism $f\colon X \ra S$ of relative dimension $1$ such that $\sZ$ is $S$-finite and $D$ is $S$-(finite \'{e}tale).  We base change $f$ along the map $\Spec(R) \ra S$ to obtain
\bitem
\m
a smooth affine $R$-scheme $C$ of pure relative dimension $1$ (base change of $X$);


\m
an $R$-finite closed subscheme $Z \subset C$ (base change of $\sZ$);

\m
an $R$-(finite \'{e}tale) closed subscheme $Y \subset Z$ (base change of $D$); 

\m
a section $s \in C(R)$ (induced by the ``diagonal'' section) such that $s|_{R[\f1r]}$ factors through $C \setminus Y$;

\m
a reductive $C$-group $\sG$ with $s^*(\sG) \cong \cG$ (base change of $\cG$) such that $\sG_{Y}$ is totally isotropic; 

\m
a maximal $C$-torus $\sT \subset \sG$ (base change of $T$) with $s^*(\sT) \cong T$; and

\m
a $\sG$-torsor $\sE$ over $C \setminus Y$ (base change of $E$) that is trivial over $C\setminus Z$ such that
\[
\tst \qq (s|_{R[\f1r]})^*(\sE) \cong E \qx{as $\cG$-torsors over $R[\f1r]$.}
\]
\eitem
We replace $Z$ by $Z \cup s$ if needed to arrange that $s \in Z(R)$. By \Cref{lem:equate-groups} (with \S\ref{conv} for the isotriviality aspect) and spreading out, there is a finite \'{e}tale cover $\wt{C}$ of some affine open neighborhood of $Z$ in $C$ such that $s$ lifts to some $\wt{s} \in \wt{C}(R)$ and $\sG_{\wt{C}} \simeq \cG_{\wt{C}}$, compatibly with an already fixed such isomorphism after pullback along $\wt{s}$. Thus, we may replace $C$ and $s$ by $\wt{C}$ and $\wt{s}$ and replace $Z$, $Y$, $\sG$, $\sE$ by their corresponding base changes to reduce to when $\sG$ is $\cG_C$. 
In this case, however, by \Cref{thm:A1}, the $\cG$-torsor $\sE$ extends to a $\cG$-torsor defined over all of $C$.  Thus, by pulling back along $s$, our $\cG$-torsor $E$ extends to a desired $\cG$-torsor $\cE$ over $R$. 
\QED
\epp

\addtocounter{footnote}{1}
\renewcommand{\thefootnote}{\arabic{footnote}}


\section{The generalized Bass--Quillen conjecture over a field}

The proof of \Cref{thm:BQ-ann} will use the following general form of Quillen~patching.

\blem[Gabber, see \cite{torsors-regular}*{Corollary 5.1.5 (b)}] \label{thm:Quillen-patching}
For a ring $A$ and a locally finitely presented $A$-group algebraic space $G$, a $G$-torsor \up{for the fppf topology} over $\bA^d_A$ descends to a $G$-torsor over $A$ if and only if it does so Zariski locally on $\Spec(A)$. \QED
\elem

\bpp[Proof of \Cref{thm:BQ-ann}] \label{pp:BQ-pf}
We have a regular ring $R$ containing a field, a totally isotropic 
reductive $R$-group $G$, and a generically trivial $G$-torsor $E$ over $\bA^d_R$. We need to show that $E$ descends to a $G$-torsor over $R$. For this, by induction on $d$, we may assume that $d = 1$. By Quillen patching of \Cref{thm:Quillen-patching}, we may assume that $R$ is local. In this key local case, we will show that $E$ is trivial.

For this, by \Cref{conj:Horrocks}, it suffices to show that $E$ is trivial on $\bA^1_R \setminus Z$ for some $R$-finite closed subscheme $Z \subset \bA^1_R$. By a limit argument, it therefore suffices to show that $E$ becomes trivial over the localization of $R[t]$ obtained by inverting all the monic polynomials. By the change of variables $x \ce t\i$, this localization is the localization of $\bP^1_R$ along the section $\infty$, and hence is isomorphic to 
\[
\tst (R[x]_{1 + xR[x]})[\f1x].
\] 
The ring $R' \ce R[x]_{1 + xR[x]}$ is regular, local, and shares its fraction field with $\bA^1_R$. In particular, the base change of $E$ to $R'$ is generically trivial. Thus, since $x$ is a regular parameter of $R'$, \Cref{thm:main-DVR}~\ref{m:mDVR-1} implies that this base change of $E$ is trivial, as desired.
 \QED
\epp













\end{document}